\documentclass[12pt]{amsart}
\usepackage[latin1]{inputenc}
\usepackage{amssymb,amsmath}

\newtheorem{thm}{Theorem}[section]
\newtheorem{cor}[thm]{Corollary}
\newtheorem{lem}[thm]{Lemma}
\newtheorem{prop}[thm]{Proposition}

\theoremstyle{definition}
\newtheorem{defn}[thm]{Definition}
\newtheorem{question}{Question}
\theoremstyle{remark}
\newtheorem{rem}[thm]{Remark}
\numberwithin{equation}{section}

\newcommand{\To}{\longrightarrow}

\newcommand{\nat}{\mathbb{N}}

\begin{document}

\baselineskip=17pt

\title{Fiber orders and compact spaces of uncountable weight}
\author{Antonio Avil\'es and Ond\v{r}ej F.K. Kalenda}
\address{University of Murcia, Facultad de Matemáticas, Campus de Espinardo, Murcia, Spain} \email{avileslo@um.es}
\address{Charles University, Faculty of Mathematics and Physics, Department of Mathematical Analysis, Sokolovsk\'a 83, 186~75 Praha 8, Czech Republic}
\email{kalenda@karlin.mff.cuni.cz}

\thanks{The first author was supported by a Marie Curie
Intra-European Felloship MCEIF-CT2006-038768 (E.U.) and research
projects MTM2005-08379 (MEC and FEDER) and S\'{e}neca 00690/PI/04}
\thanks{The second author was supported by Research project MSM 0021620839 financed by MSMT
and by Research grant GA\v{C}R 201/06/0018}



\begin{abstract}
We study an order relation on the fibers of a continuous map and
its application to the study of the structure of compact spaces of
uncountable weight.
\end{abstract}

\maketitle

\section{Introduction and main results}

This work is motivated by the following general problem: Given two
compact convex sets $K$ and $L$ (sitting in some locally convex
linear topological spaces), are $K$ and $L$ homeomorphic? When $K$
and $L$ are metrizable (that is, they have countable weight) the
well known Keller's theorem, cf. \cite{BesPelInfDim}, implies that
$K$ and $L$ are homeomorphic if and only if they have the same
dimension. Thus, when restricting our attention to compact sets of
countable weight, only one topological invariant has to be
computed to
answer our question: the dimension, ranging from 0 to $\omega$.

When we pass to the case when the weight is uncountable, the
situation is not that simple. A number of usual topological
invariants, like chain conditions, cardinal functions,
functional-analytic properties, etc. can be used to identify many
different types of compact convex sets. Just to recall an
elementary example, we may compare an uncountable product of
intervals $[0,1]^\kappa\subset\mathbb{R}^\kappa$ with the ball
$B(\kappa)$ of the Hilbert space $\ell_2(\kappa)$ in the weak
topology. In $B(\kappa)$ we may find an uncountable family of
disjoint open sets but $[0,1]^\kappa$ has the countable chain
condition. Another argument would be that $B(\kappa)$ cannot be
homeomorphic to an uncountable product since it contains
$G_\delta$ points (we will obtain in this paper a much subtler
fact: $B(\kappa)$ is not homeomorphic even to a finite product of
compact spaces
of uncountable weight).

Very often, however, the standard topological technology is not so
helpful as in the mentioned case of $B(\kappa)$ and
$[0,1]^\kappa$. An example of this is when we restrict our
attention to weakly compact sets of the Hilbert space
$\ell_2(\kappa)$. The first example a nonmetrizable weakly compact convex set not homeomorphic to $B(\kappa)$ may be traced back to constructions of Corson and Lindenstrauss \cite{CorsonLindenstrauss,Lindenstrauss}, who provided such a set in which all points are $G_\delta$. Such sets, however, cannot be symmetric and to the best of our knowledge, only recently the first
author~\cite{AvilesHilbertBall} provided a first example of an absolutely
convex weakly compact subset of $\ell_2(\kappa)$ of weight
$\kappa$ which is not homeomorphic to $B(\kappa)$. This was done
by proving that $B(\kappa)$ satisfies a certain chain condition of
Ramsey type introduced by Bell~\cite{BellRamsey} and constructing
ad hoc a compact convex set failing such property.

Let us provide now some natural examples of compact convex sets,
all of them indeed representable as weakly compact convex subsets
of $\ell_2(\kappa)$, for which apparently the standard techniques
from topology give us no clue about the problem whether they are
homeomorphic to each other or not. Letter $\kappa$ always denotes an uncountable cardinal.

\begin{itemize}
\item The ball of the Hilbert space $B(\kappa) = \{
x\in\ell_2(\kappa) : \|x\|_2\leq 1\}$.\smallskip

\item The space $P(A(\kappa))$ of Radon probability measures on
$A(\kappa)=\kappa\cup\{\infty\}$, the one-point compactification of the discrete set $\kappa$.
\smallskip

\item The spaces $P(A(\kappa)^n)$, $2\leq n<\omega$.
\smallskip
\item The spaces $P(\sigma_n(\kappa))$ of probability measures on
$\sigma_{n}(\kappa) = \{x\in\{0,1\}^\kappa : |\mbox{supp}(x)|\leq n\}$, $2\leq n <\omega$.
\smallskip
\item The finite and countable powers of the previous spaces.

\end{itemize}

We shall develop some new tools which will allow us to conclude
that all these spaces are not homeomorphic to each other, with
perhaps the exception of $B(\kappa)$ and $P(A(\kappa))$ for which
our techniques are unable to determine whether they are
homeomorphic or not. We also studied other examples, not
embeddable into a Hilbert space, namely the compact sets
$P([0,\kappa]^n)^m$ for $\kappa$ uncountable regular cardinal and
$n,m\in\mathbb{N}$. In addition, we will obtain other applications
concerning
the structure of these spaces, regarding the two following kind of questions:

\begin{itemize}

\item The classification of the points of a compact space $K$, that is, for which points $x,y\in K$ there exists a homeomorphism $f:K\To K$ such that $f(x)=y$.

\item When a compact space $K$ can be homeomorphic to some power compact of the form $L^n$, or when it can be homeomorphic to a product of the form $L_1\times\cdots\times L_n$.

\end{itemize}

The way to address all these questions goes through the beautiful
technique of Shchepin of inverse limits and the spectral theorem
developped in \cite{Shchepin1} and \cite{Shchepin2}. We explain
this in detail in Section~\ref{sectionspectral}, but roughly
speaking, given a compact space $K$ of uncountable weight, this
technique allows to study the topological structure of $K$ by
studying the continuous surjections $p:X\To Y$ for $X$ and $Y$
quotients of $K$ of countable weight. And here comes the key idea
of our work, to study a certain
preorder relation induced on the fibers of a continuous map:

\begin{defn} Let $f:K\To L$ be a continuous map and $x\in L$. We
define a preorder relation $\leq$ on the fiber $f^{-1}(x)$ by
letting $s\leq t$ if and only if for every neighborhood $U$ of $s$
there exists a neighborhood $V$ of $t$ such that $f(V)\subset
f(U)$.
\end{defn}

In other words, $s\leq t$ if and only if
$$\{f(U) : U\text{ is a neighborhood of }s\}\subset \{f(V) : V\text{ is a neighborhood of
}t\},$$ if and only if $f^{-1}f(U)$ is a neighborhood of $t$ for
every neighborhood $U$ of $s$. We shall call $\mathbb{F}_x(f) =
f^{-1}(x)$ to the fiber of $x$ endowed with the preorder $\leq$
(and also with its inherent topology, though we shall not use the
topological structure here). We denote by $\mathbb{O}_x(f) =
\mathbb{F}_x(f)/\sim$ the ordered set obtained by making a
quotient by the equivalence relation $t\sim s$ $\iff$ $t\leq s$
and $s\leq
t$.

 In his study of the spaces $exp_n(2^\kappa)$
\cite{Shchepin1}, Shchepin considered what in our language would
be the cardinality of $\mathbb{O}_x(f)$. This was already useful
in that \emph{discrete} context but not in spaces like convex
sets, where one needs to consider the ordered structure of
$\mathbb{O}_x(f)$ to get some information.

Let us indicate how fiber orders may be helpful in the problem of
classification of points of a compact space, and in the
homeomorphic classification of compact sets. Consider a compact
space $K$ of uncountable weight and a point $x\in K$. We can
consider then the family of all fiber orders of type
$\mathbb{O}_{p_L(x)}(q)$ for every continuous surjection $q:L'\To
L$ between metrizable quotients of $K$ with projections $p_{L'}:K\To
L'$, $p_L:K\To L$, $qp_{L'}=p_L$. This collection of ordered sets may
be in principle rather complicated, but in the examples that we
deal with it happens that \emph{almost} all these sets are order
isomorphic to the same ordered set that we can call
$\mathbb{O}_x(K)$. For instance, for a finite power of the ball of
the nonseparable Hilbert space
$B(\kappa)$ we get the following picture:

\begin{thm}\label{OxBk}
Let $K=B(\kappa)^n$ and $x=(x_1,\ldots,x_n)\in K$. Let $r =|\{i :
\|x_i\| < 1\}|$. Then $\mathbb{O}_x(K)\cong [0,1]^r$.
\end{thm}

We view $[0,1]^r$ as an ordered set endowed with the pointwise
order, i.e.
$$(t_1,\ldots,t_r)\leq (s_1,\ldots,s_r)\mbox{ iff }t_i\leq s_i
\mbox{ for every }i.$$
 Notice other consequences of this result other than
the fact that the finite powers of the ball are nonhomeomorphic.
It is a standard fact that the points of $B(\kappa)^n$ whose all
coordinates belong to the sphere are the $G_\delta$ points of
$B(\kappa)^n$ and hence, topologically different from the rest. We
obtained something much less evident: that points with different
number of coordinates in the sphere are topologically different.
This is a complete classification of the points of $B(\kappa)^n$
because if two points have the same number of coordinates in the
sphere, then there is an automorphism of $B(\kappa)^n$ which moves
one to the other.

Apart from the euclidean ball, the other spaces that we studied
are spaces of probability measures on scattered spaces. We
developped a general method for computing fiber orders in these
cases, which constitutes the part of our work which is technically
the most involved. One of the key steps in this task is our
Lemma~\ref{imageofneighborhood} which probably has an independent
interest. Every Radon probability measure on a scattered compact
space is discrete, thus a certain (finite or infinite) convex
combination of Dirac measures $\delta_x$. The following result
(which follows immediately from Theorem~\ref{sumofprobabilities} below)
reduces the computation of fiber orders in spaces
$P(K)$ to Dirac measures:

\begin{thm}\label{Oxsumadiraku}
Let $K$ be a scattered compact and $\mu=\sum_{i\in I} r_i\delta_{x_i}\in P(K)$, where $x_i\in K$ are pairwise distinct and $r_i>0$ for $i\in I$. Then $$\mathbb{O}_\mu(P(K))\cong\prod_{i\in I} \mathbb{O}_{\delta_{x_i}}(P(K)).$$
\end{thm}

The picture of the fiber orders of Dirac measures in our examples
of probability measures spaces is the following:

\begin{thm}\label{symmetric}
Let $K=\sigma_n(\kappa)$ and $x\in K$. Set $k=n-|x|$. Then
$\mathbb{O}_{\delta_x}(P(K))\cong \{(t_1,\ldots,t_k)\in [0,1]^k :
t_1\leq\ldots\leq t_k\}$ where the order is defined as $t\leq s$
if and only if
$t_j\leq s_j$ for every $j$.
\end{thm}

In the next result, we denote by $2^k$ the power set of
$\{1,\ldots,k\}$.

\begin{thm}\label{polyhedric} Let $K=A(\kappa)^n$ or $K=[0,\omega_1]^n$,
and let $x=(x_1,\ldots,x_n)\in K^n$. Set $k$ to be the number of
coordinates of $x$ which are not $G_\delta$-points of $K$. Then
$\mathbb{O}_{\delta_x}(P(K))\cong\{ (t_A)_{A\in 2^k}\in
[0,1]^{2^k} : \sum_{A\in 2^k}t_A =1\}$ endowed with the order
$(t_A)\leq (s_A)$ if and only if $\sum_{A\in\mathcal{A}}t_A \leq
\sum_{A\in \mathcal A}s_A$ for every upwards closed family
$\mathcal A$ of subsets of $\{1,\ldots,k\}$.
\end{thm}

A similar statement as Theorem~\ref{polyhedric} holds for compact
spaces $K=[0,\tau]$ with $\tau$ an uncountable regular cardinal,
but for a modified version of the ordered sets $\mathbb{O}_x(L)$
relative to the cardinal $\tau$. Finally we state the kind of
results that we prove using these techniques that refer to
decomposition of
compact spaces as products:

\begin{thm}\label{powerBall}
Let $K=B(\kappa)$, $P(\sigma_n(\kappa))$, $P(A(\kappa)^n)$ or
$P([0,\tau]^n)$ for $\tau$ an uncountable regular cardinal, and
let $k,m\in\mathbb{N}$. Suppose that there exists a compact $L$
such that $K^k \approx L^m$. Then, $k$ is a multiple of $m$.
\end{thm}

\begin{thm}\label{productball}
Let $K=B(\kappa)$, $P(A(\kappa))$ or $P([0,\omega_1])$ and let
$n,m\in\mathbb{N}$. Suppose that $L_1,\ldots,L_m$ are compact
spaces of uncountable weight such that $K^n \approx \prod_{i=1}^m
L_i$. Then,
$m\leq n$.
\end{thm}

\begin{thm}\label{producttau}
Let $\tau$ be a regular cardinal, $n,m$ natural numbers, and
$L_1,\ldots,L_m$ compact spaces of weight $\tau$. If
$P([0,\tau])^n\approx L_1\times\cdots\times L_m$, then $m\leq n$.
\end{thm}

We make two remarks about these results. First,
our methods do not allow to decide whether these compact spaces can be expressed as a nontrivial product with one factor metrizable. This appears not to be an easy question.
Using a result of \cite{schori} and its variants, the second author \cite{Kalendaproducts} has obtained that $P(K)$ is homeomorphic to $P(K)\times[0,1]$, for any compact scattered $K$. However it is unknown to us whether $B(\kappa)$ is homeomorphic to $B(\kappa)\times [0,1]$. Second, the first author \cite{Avilesproducts} has obtained with different techniques an improvement of Theorem~\ref{productball}: If $B(\kappa)^n$ maps continuously onto a product of nonmetrizable compacta of the form $\prod_{i=1}^m L_i$, then $m\leq n$. These techniques do not apply to the case of Theorem~\ref{powerBall} for $K\neq B(\kappa)$, and actually $P(A(\kappa)^n)$ and $P(\sigma_n(\kappa))$ map continuosly onto $B(\kappa)^n$.

\section{Spectral theory}\label{sectionspectral}

In this section, we summarize in a self-contained way what we need
about spectral theory, which is essentially taken from
\cite{Shchepin1} and \cite{Shchepin2}. We also introduce the
invariants $\mathbb{F}_x(K)$ and $\mathbb{O}_x(K)$, which play a
central
role in the paper.

Let $K$ be a compact space. We denote by $\mathcal{Q}(K)$ the set
all Hausdorff quotient spaces of $K$, that is the set all
Hausdorff compact spaces of the form $K/E$ endowed with the
quotient topology, for $E$ an equivalence relation on $K$. An
element of $\mathcal{Q}(K)$ can be represented either by the
equivalence relation $E$ or by the quotient space $L=K/E$ together
with the canonical projection $p_L:K\To L$.

On the set $\mathcal{Q}(K)$ there is a natural order relation. In
terms of equivalence relations $E\leq E'$ if and only if
$E'\subset E$. Equivalently, in terms of the quotient spaces,
$L\leq L'$ if and only if there is a continuous surjection
$q:L'\To L$ such that $q p_{L'} = p_L$. The set $\mathcal{Q}(K)$
endowed with this order relation is a complete semilattice, that
is, every subset has a least upper bound or supremum: if
$\mathcal{E}$ is a family of equivalence relations of
$\mathcal{Q}(K)$, its least upper bound is the relation given by
$x E_0 y$ if and only if $x E y$ for all $E\in\mathcal{E}$, in other words
$E_0 =\sup\mathcal{E} = \bigcap \mathcal{E}$. It is easy to check
that $E_0$ gives a
Hausdorff quotient if each element of $\mathcal{E}$ does.

 Let
$\mathcal{Q}_\omega(K)\subset \mathcal{Q}(K)$ be the family of all
quotients of $K$ which have countable weight. Notice that
$\sup\mathcal{A}\in \mathcal{Q}_\omega(K)$ for every countable
subset $\mathcal{A}\subset\mathcal{Q}_\omega(K)$ and also that
$\sup \mathcal{Q}_\omega(K)=K$. A family $\mathcal{S}\subset
\mathcal{Q}_\omega(K)$ is called cofinal if for every $L\in
\mathcal{Q}_\omega(K)$ there exists $L'\in\mathcal{S}$ such that
$L\leq L'$. The family $\mathcal{S}$ is called a
$\sigma$-semilattice if for every countable subset
$\mathcal{A}\subset\mathcal{S}$, the least upper bound of
$\mathcal{A}$ belongs to $\mathcal{S}$.

\begin{thm}[A version of Shchepin's spectral theorem]\label{ourspectraltheorem}
Let $K$ be a compact space of uncountable weight and let
$\mathcal{S}$ and $\mathcal{S'}$ two cofinal $\sigma$-semilattices
in $\mathcal{Q}_\omega(K)$. Then $\mathcal{S}\cap\mathcal{S'}$ is
also a cofinal $\sigma$-semilattice in $\mathcal{Q}_\omega(K)$.
\end{thm}

Proof: The point is in proving that $\mathcal{S}\cap\mathcal{S'}$
is cofinal. Let $L_0\in \mathcal{Q}_\omega(K)$ be arbitrary. Since
$\mathcal{S}$ is cofinal there exists $L_1\in\mathcal{S}$ with
$L_0\leq L_1$, similarly find $L_2\in \mathcal{S}'$ with $L_1\leq
L_2$, and continue by induction an increasing sequence with
$L_{2n+1}\in\mathcal{S}$, $L_{2n}\in \mathcal{S}'$. Finally
$L=\sup\{L_n : n<\omega\}\in\mathcal{S}\cap\mathcal{S}'$ since
both sets are $\sigma$-semilattices.$\qed$

It is not so obvious to check whether a given $\sigma$-semilattice
is cofinal, so this theorem must be applied together with the
following criterion:

\begin{lem}\label{spectralfactor}
Let $K$ be a compact space of uncountable weight and $\mathcal{S}$
a $\sigma$-semilattice in $\mathcal{Q}_\omega(K)$. Then,
$\mathcal{S}$
is cofinal if and only if $\sup\mathcal{S} = K$.
\end{lem}

Proof: If $\mathcal{S}$ is cofinal, then $\sup\mathcal{S} =
\sup\mathcal{Q}_\omega(K) = K$. Conversely, suppose that
$\sup\mathcal{S}=K$. Consider the family $\mathcal{A}$ of all
continuous functions $f:K\To\mathbb{R}$ such that there exists
$L\in\mathcal{S}$ such that $f$ factors through $p_L:K\To L$, that
is, there exists $\hat{f}:L\To\mathbb{R}$ with $f=\hat{f}p_L$. As
$\mathcal{S}$ is a $\sigma$-semilattice, $\mathcal{A}$ is a
subalgebra of the algebra $C(K)$ of real-valued continuous
functions on $K$. Clearly, constant functions belong to
$\mathcal{A}$ and since $\sup\mathcal{S}=K$, $\mathcal{A}$
separates the points of $K$. Hence, by the Stone-Weierstrass
theorem every $f\in C(K)$ is the limit of a sequence of functions
from $\mathcal{A}$. But indeed $\mathcal{A}$ is closed under
limits of sequences, namely if $f_n$ factors throught
$L_n\in\mathcal{S}$, then $\lim f_n$ factors through $\sup\{L_n :
n<\omega\}\in\mathcal{S}$. We conclude that $\mathcal{A}=C(K)$.
Now, if $p:K\To L$ is an arbitrary element of
$\mathcal{Q}_\omega(K)$, then we can take an embedding
$L\subset\mathbb{R}^\omega$ and consider the functions
$e_np:K\To\mathbb{R}$ obtained by composing with the coordinate
functions $e_n:\mathbb{R}^\omega\To\mathbb{R}$. For every $n$ we
know, since $\mathcal{A}=C(K)$, that there exists
$L_n\in\mathcal{S}$ such that $e_n p$ factors through $L_n$.
Finally, this implies that $p$ factors through $L_\infty =
\sup\{L_n : n<\omega\}$, so $L\leq
L_\infty\in\mathcal{S}$.$\qed$

The importance of this machinery is that it allows to study a
compact space of uncountable weight through the study of a cofinal
$\sigma$-semilattice of metrizable quotients, and particularly
through the natural projections between elements of the
$\sigma$-semilattice. In this way, the study of compact spaces of
uncountable weight is related to the study of continuous
surjections between compact spaces of countable weight. The
following language will be
useful:

\begin{defn}
Let $K$ be a compact space of uncountable weight and let
$\mathcal{P}$ be a property. We say that the $\sigma$-typical
surjection of $K$ satisfies property $\mathcal P$ if there exists
a cofinal $\sigma$-semilattice $\mathcal{S}\subset
\mathcal{Q}_\omega(K)$ such that for every $L< L'$ elements of
$\mathcal{S}$, the natural projection $p:L'\To L$ satisfies
property $\mathcal{P}$.
\end{defn}

The consequence of the spectral theorem is that the fact whether
the $\sigma$-typical surjection of $K$ has a certain property can
be checked on any given cofinal $\sigma$-semilattice, namely:

\begin{thm}
Let $K$ be a compact space of uncountable weight, let
$\mathcal{P}$ be a property, and let $\mathcal{S}$ be a fixed
cofinal $\sigma$-semilattice in $\mathcal{Q}_\omega(K)$. Then, the
$\sigma$-typical surjection of $K$ has property $\mathcal{P}$ if
and only if there exists a cofinal $\sigma$-semilattice
$\mathcal{S}'\subset\mathcal{S}$ such that for every $L< L'$
elements of $\mathcal{S}'$, the natural projection $p:L'\To L$
satisfies property $\mathcal{P}$.
\end{thm}

The main kind of properties $\mathcal P$ that we shall be
interested concern the fiber orders of the surjections and the
order relation that we defined on them. Given a point $x\in K$, we
can study properties of the point $x$ by looking to fiber order of
$p_L(x)$ in the $\sigma$-typical $p:L'\To L$. It may be a useful
language to call $\mathbb{F}_x(K)$ to this $\sigma$-typical fiber,
which we certainly cannot define as a concrete set, but rather as
an abstract object of which we can predicate some properties.

\begin{defn}
Let $K$ be a compact space of uncountable weight and $x\in K$ and
let $\mathcal{P}$ be a property. We say that $\mathbb{F}_x(K)$ has
property $\mathcal P$ if $\mathbb{F}_{p_L(x)}(p)$ has
property $\mathcal P$ for the $\sigma$-typical surjection $p:L'\To L$.
\end{defn}

In a similar way we shall talk about $\mathbb{O}_x(K)$. It is
worth to notice that a point $x$ is a $G_\delta$-point of $K$ if
and only if $|\mathbb{F}_x(K)|=1$. In other words, the information
given by $\mathbb{F}_x(K)$ is trivial only when $x$ is a
$G_\delta$ point of $K$. Namely, if $x$ is a $G_\delta$-point of
$K$ then there is a continuous function $f:K\To [0,1]$ such that
$x=f^{-1}(0)$. Then, $f$ can be viewed as an element $L_0\in
\mathcal{Q}_\omega(K)$ and we find that
$|\mathbb{F}_{p_L(x)}(p)|=1$ for all $L'>L>L_0$, $p:L'\To L$.
Conversely, if $|\mathbb{F}_x(K)|=1$ then we can find
$L\in\mathcal{Q}_\omega(K)$ such that $x=p_L^{-1}(p_L(x))$.
Another elementary example is the compact $K=L^\kappa$ where $L$
is a metrizable compact. In this case, one can see as an excercise
that for every $x\in K$, $\mathbb{F}_x(K)$
is homeomorphic to $L^\omega$ and $|\mathbb{O}_x(K)|=1$.

\section{Decomposition into products}

In this section, apart from providing some basic facts that will
be needed in the sequel, we prove two results,
Theorem~\ref{Knoroots} and Theorem~\ref{Knodecomposition}, which
establish some sufficient conditions in terms of fiber orders in
order that a compact $K$ cannot be decomposed as product of other
spaces in a certain way. In further sections, when computing the
fiber orders of specific spaces, we will find that several compact
spaces satisfy the assumptions of these results.

\begin{defn}
Let $P$ be a set and $\leq$ be a binary relation on $P$. We say
that $(P,\leq)$ is a preordered set if \begin{enumerate} \item
$t\leq t$ for every $t\in P$, \item If $t\leq s$ and $s\leq u$,
then $t\leq u$, for every $t,s,u\in P$.
\end{enumerate}
If, moreover, we have that for every $t,s\in P$, if $t\leq s$ and
$s\leq t$ then $t=s$, then we say that $(P,\leq)$ is an ordered
set. An ordered set $(O,\leq)$ is said to be linearly ordered if
for every $t,s\in O$, either $t\leq s$ or $s\leq t$.
\end{defn}

There is a canonical way of constructing an ordered set from a
given preordered set $(P,\leq)$, namely we consider the
equivalence relation on $P$ given by $t\sim s$ iff $t\leq s$ and
$s\leq t$, and then the quotient set $P/\sim$ is an ordered set
when endowed with the relation induced from $P$. We call this the
ordered set associated to $P$. When we write $p<q$ in a preordered
set, it means that $p\leq q$ but $q\not\leq p$.

An isomorphism between the preordered sets $P$ and $Q$ is a
bijection $f:P\To Q$ such that $f(t)\leq f(s)$ if and only if
$t\leq s$.

\begin{defn}
Let $\{Q_i: i\in I\}$ be a family of preordered sets. The product
of this family is the preordered set whose underlying set is the
cartesian product $\prod_{i\in I}Q_i$ endowed by the preorder
relation given by: $(t_i)_{i\in I}\leq
(s_i)_{i\in I}$ if and only if $t_i\leq s_i$ for every $i\in I$.
\end{defn}

The product of an empty family of preordered sets is considered to
be a singleton, with its only possible preordered structure. The
product operation of preordered sets arises naturally in the
context of fiber orders at least in two different situations,
related to probability measures (cf.
Theorem~\ref{sumofprobabilities}) and to products of compact
spaces:

\begin{prop}\label{productbasics} Let $\{f_i:K_i\To L_i : i\in I\}$ be a family of
continuous surjections, let $f:\prod_{i\in I} K_i\To \prod_{i\in I} L_i$ be its
product and let $x=(x_i)_{i\in I}$ be a point of $\prod_{i\in I} L_i$.
Then, the natural map $\mathbb{F}_x(f)\To \prod_{i\in
I}\mathbb{F}_{x_i}(f_i)$ is an order-isomorphism. In particular,
$\mathbb{F}_x(f)\cong \prod_{i\in I}\mathbb{F}_{x_i}(f_i)$ and
$\mathbb{O}_x(f)\cong \prod_{i\in
I}\mathbb{O}_{x_i}(f_i)$.
\end{prop}

The proof of this statement is straightforward. If we have $K$ a
finite or countable product of compact spaces, then a cofinal
$\sigma$-semilattice in $\mathcal{Q}_\omega(K)$ is formed by all
quotients of countable weight of $K$ which can be expressed as the
product of a quotient of every factor. In this way, we see that
the fibers of the $\sigma$-typical surjection of the product are
the product of the fibers of the $\sigma$-typical surjection of
every factor. We are thus allowed to write expressions like for
instance $\mathbb{F}_{(x,y)}(K\times L) \cong
\mathbb{F}_x(K)\times \mathbb{F}_y(L)$ or
$\mathbb{F}_{(x_1,x_2,\ldots)}(\prod_{n<\omega}K_n)\cong \prod_{n<\omega}\mathbb{F}_{x_n}(K_n)$.

\begin{defn}
An ordered set $O$ is called \emph{irreducible} if whenever $O$ is
isomorphic to a product $Q\times R$ we have that either $Q$ or $R$
is a singleton.
\end{defn}

An elementary example of an irreducible ordered set is a linearly
ordered set. An ordered set $O$ is called \emph{connected} if
whenever it is expressed as the disjoint union of two nonempty
subsets $O=A\cup B$, there exists $a\in A$ and $b\in B$ such that
either $a\leq b$ or $b\leq a$. All the ordered sets that appear in
this note happen to be connected since indeed they have a minimum. The following Theorem~\ref{productofordered} and its
Corollary~\ref{productofirreducibles} are due to
Hashimoto~\cite{Hashimoto} and assert that any two decompositions
of a connected ordered set as product have a common refinement,
and consequently, a decomoposition of a connected ordered set as a
product of irreducible ordered sets is unique. Among other
applications, this is a useful criterion to decide immediately
that two given ordered sets are not isomorphic.

\begin{thm}\label{productofordered}
Let $\mathcal{O}$ be a connected ordered set, $\{O_i : i\in I\}$
and $\{Q_j : j\in J\}$ two families of ordered sets such that
$O\cong\prod_{i\in I}O_i\cong \prod_{j\in J}Q_j$. Then, there is a
further family $\{Z_{ij}: (i,j)\in I\times J\}$ such that $O_i
\cong\prod_{j\in J}Z_{ij}$ for every $i\in I$, and
$Q_j\cong\prod_{i\in I}Z_{ij}$ for every $j\in J$.
\end{thm}

\begin{cor}\label{productofirreducibles}
Let $\mathcal{O}$ be a connected ordered set, $\{O_i : i\in I\}$ a
family of irreducible ordered sets and $\{Q_j : j\in J\}$ a family
of arbitrary ordered sets. Assume that $O\cong\prod_{i\in
I}O_i\cong \prod_{j\in J}Q_j$. Then, there is a partition
$I=\bigcup_{j\in J}F_j$ of the set $I$
such that $Q_j\cong \prod_{i\in I_j}O_i$ for every $j\in J$.
\end{cor}

In the sequel we shall make use of the following terminology: Two
continuous maps $f:U\To V$ and $f':U'\To V'$ are said to be
homeomorphic if there exists homeomorphisms $u:U\To U'$ and
$v:V\To V'$ such that $vf = f'u$.

\begin{thm}\label{Knoroots}
Let $K$ be a compact space of uncountable weight and let $O$ be a
connected irreducible ordered set. Assume that there is $x\in K$
such that $\mathbb{O}_x(K)\cong O$ and that
$\mathbb{O}_y(K)\not\cong O^k$ for each $y\in K$ and each $k>1$.
If $K^n\approx L^m$ for some natural numbers $n$, $m$ and some
space $L$, then $n$ is a multiple of $m$.
\end{thm}

\begin{rem} Note that the assertion $\mathbb{O}_y(K)\not\cong O^k$ is not the negation of
$\mathbb{O}_y(K)\cong O^k$. It rather means that for the
$\sigma$-typical surjection $p$ we have
$\mathbb{O}_{y}(p)\not\cong O^k$ (another remark about notation:
we write $\mathbb{O}_{y}(p) = \mathbb{O}_{y'}(p)$, where $y'$ is
the projection of $y$ on the range of $p$).
\end{rem}

Proof of Theorem~\ref{Knoroots}: Along this proof, it is important
to have in mind that if $\mathcal{S}$ is a cofinal
$\sigma$-semilattice in $\mathcal{Q}_\omega(X)$, then
$\mathcal{S}^k = \{p_Z^k:X^k\To Z^k: Z\in\mathcal{S}\}$ is a
cofinal $\sigma$-semilattice in $\mathcal{Q}_\omega(X^k)$,
$k\leq\omega$. Assume that $n$ is not a multiple of $m$ and that
$K^n\approx L^m$. Choose $x\in K$ with $\mathbb{O}_x(K)\cong O$.
By Proposition~\ref{productbasics} we get
$\mathbb{O}_{(x,\dots,x)}(K^n)\cong O^n$. Let
$w=(w_1,\ldots,w_m)\in L^m$ be the point corresponding to
$(x,\dots,x)$ by the homeomorphism. Then, of course,
$\mathbb{O}_w(L^m)\cong O^n$. Further, by
Proposition~\ref{productbasics} we have that for the
$\sigma$-typical surjection $q$ of $L$,
$$O^n\cong\mathbb{O}_{w}(q^m)\cong\prod_{i=1}^m\mathbb{O}_{w_i}(q).$$
Using Corollary~\ref{productofirreducibles} and the fact that $n$
is not a multiple of $m$, we get that for the $\sigma$-typical
surjection $q$ of $L$ there is $k\in\{1,\dots,m\}$ and $n/m<s\leq
n$ such that $\mathbb{O}_{w_k}(q)\cong O^s$. It follows that there
are $k\in\{1,\dots,m\}$ and $n/m<s\leq n$ such that in each
cofinal $\sigma$-semilattice in $L$ there is some surjection $q$
with $\mathbb{O}_{w_k}(q)\cong O^s$ (this follows from
Theorem~\ref{ourspectraltheorem}: if not, for each $k,s$ there
would be the corresponding cofinal $\sigma$-semilattice
$\mathcal{S}_{k,s}$, and then $\bigcap\mathcal{S}_{k,s}$ gives a
contradiction). Set $\tilde w=(w_k,\dots,w_k)$ and let
$y=(y_1,\dots,y_n)\in K^n$ correspond by the homeomorphism to
$\tilde w$. By our assumptions there is a cofinal
$\sigma$-semilattice $\mathcal{T}\subset\mathcal{Q}_\omega(K)$
such that $\mathbb{O}_{p(y_i)}(p)\not\cong O^j$ for any surjection
$p$ inside $\mathcal{T}$, $i=1,\dots,n$ and $j>1$. Consider the
cofinal $\sigma$-lattice $\mathcal{U}=\{Z\in \mathcal{Q}_\omega(L)
: Z^m\in \mathcal{T}^n\}$. By the previous argument, there exists
a surjection $q$ inside $\mathcal{U}$ such that
$\mathbb{O}_{w_k}(q)\cong O^s$. The surjection $q^m$ corresponds
to a surjection $p^n$ inside $\mathcal{T}^n$ for which we have
that:
$$\prod_{i=1}^n\mathbb{O}_{y_i}(p)\cong\mathbb{O}_{y}(p^n)\cong O^{sm}$$
As $sm>n$, by Corollary~\ref{productofirreducibles} we get that
$\mathbb{O}_{p(y_i)}(p)\cong O^j$ for some $i\in\{1,\dots,n\}$ and
some $j>1$, a contradiction.
$\qed$

\begin{thm}\label{Knodecomposition}
Let $K$ be a compact space of uncountable weight, $n$ a natural
number, and $O$ a connected irreducible ordered set with $|O|>1$.
Assume that the $\sigma$-typical surjection of $K$,
$p:X\To Y$ has the following properties:

\begin{enumerate}

\item For every $y\in Y$, $\mathbb{O}_y(p)$ is a connected ordered
set.

\item There is no point $y\in Y$ with $\mathbb{O}_y(p)\cong
O\times P$ with $|P|>1$.

\item There exists a point $x\in Y$ such that
$\mathbb{O}_x(p)\cong O$.

\item For any point $x\in Y$ with $\mathbb{O}_x(p)\cong O$ the
preordered set $\mathbb{F}_x(p)$ has an equivalence class which
is a
singleton.
\end{enumerate}

Then, if $L_1,\ldots,L_m$ are compact spaces of uncountable weight
such that $K^n\approx
L_1\times\cdots\times L_m$, then $m\leq n$.
\end{thm}

Proof: Let $\mathcal{S}$ be a cofinal $\sigma$-semilattice in
$\mathcal{Q}_\omega(K)$ in which all the natural projections
satisfy properties (1) to (4). Let $\mathcal{S}^n$ be the cofinal
$\sigma$-semilattice in $\mathcal{Q}_\omega(K^n)$, like defined in
the proof of Theorem~\ref{Knoroots}. Consider $\mathcal{T}$ the
cofinal $\sigma$-semilattice in
$\mathcal{Q}_\omega(L_1\times\cdots\times L_m)$ whose elements are
the quotients of $L_1\times\cdots\times L_m$ which are products of
quotients of each coordinate, that is, of the form
$$q_1\times\cdots\times q_m: L_1\times\cdots\times L_m\To
Z_1\times\cdots\times Z_m,$$ for $q_i:L_i\To Z_i$ element of
$\mathcal{Q}_\omega(L_i)$. Since $K\approx L_1\times\cdots\times
L_m$, the $\sigma$-semilattices $\mathcal{S}^n$ and $\mathcal{T}$
can be viewed as cofinal $\sigma$-semilattices of metrizable
quotients over the same compact, so by
Theorem~\ref{ourspectraltheorem} they intersect in a further
cofinal $\sigma$-semilattice, and in particular, we can find a
natural projection inside $\mathcal{S}^n$, $p^n:X^n\To Y^n$ and a
natural projection inside $\mathcal{T}$, $q=q_1\times\cdots\times
q_m:Z_1\times\cdots\times Z_m\To W_1\times\cdots\times W_m$ which
are homeomorphic. Of course, we have enough freedom to choose it
in such a way that $W_i\neq Z_i$ for every $i$. Consider a point
$w=(w_1,\ldots,w_m)$ in $W_1\times\cdots\times W_m$ which
corresponds by the homeomorphism to a point $(x,x,\ldots)\in Y^n$
with $\mathbb{O}_x(p)\cong O$. Thus,
$\prod_{r=1}^m\mathbb{O}_{w_r}(q_r)\cong O^n$. After reordering if
necessary, by Corollary~\ref{productofirreducibles} we know that
$|\mathbb{O}_{w_r}(q_r)|=1$ for $r>n$.

\emph{Claim A}: $|\mathbb{O}_{v}(q_r)|=1$ for every $v\in W_r$ and
every $r>n$.

\emph{Proof of the claim}: Suppose for instance that
there exists $v\in W_{n+1}$ with $|\mathbb{O}_{v}(q_{n+1})|>1$.
Let $w'=(w_1,\ldots,w_n,v,w_{n+2},\ldots)$ and
$y=(y_1,\ldots,y_n)\in Y^n$ which corresponds to $w'$ by the
homeomorphism. Then
$\prod_1^n\mathbb{O}_{y_i}(p)\cong\mathbb{O}_{w'}(q)\cong
O^n\times\mathbb{O}_{v}(q_{n+1})$. Using
Theorem~\ref{productofordered} and the fact that $O$ is
irreducible, we conclude that there must exist $i$ such that
$\mathbb{O}_{y_i}(p)$ is isomorphic to something of the form
$O\times P$ with $|P|>1$, which is a contradiction.

\emph{Claim B}: $|\mathbb{F}_{v}(q_r)|=1$ for every $v\in W_r$ and
every $r>n$.

\emph{Proof of the claim}: Suppose for instance that
there exists $v\in W_{n+1}$ with $|\mathbb{F}_{v}(q_{n+1})|>1$,
let $w'=(w_1,\ldots,w_n,v,w_{n+2},\ldots)$ and
$y=(y_1,\ldots,y_n)\in Y^n$ which corresponds to $w'$ by the
homeomorphism. We know by Claim A that
$|\mathbb{O}_{v}(q_{n+1})|=1$, which means that
$\mathbb{F}_{v}(q_{n+1})$ consists of one equivalence class which
is not a singleton. By Proposition~\ref{productbasics}, this
translates into the fact that $\mathbb{F}_{w'}(q)\cong
\prod_1^n\mathbb{F}_{y_i}(p)$ has no equivalence class which is
a singleton, and this further implies that for some $i$,
$\mathbb{F}_{y_i}(p)$ has no equivalence class which is  a
singleton. Moreover, $\prod_1^n\mathbb{O}_{y_i}(p)\cong
\mathbb{O}_{w'}(q)\cong O^n$, so by
Corollary~\ref{productofirreducibles} and our hypothesis~(2),
$\mathbb{O}_{y_i}(p)\cong O$ for every $i$. In this way, we found
a contradiction with our hypothesis~(4).

Finally, notice that $|\mathbb{F}_{v}(q_r)|=1$ simply means that
$q_r$ is one-to-one for $r>n$, that is $Z_r=W_r$. Since we
supposed that $Z_r\neq W_r$ for all $r$, we conclude that
$m\leq n$.$\qed$

\begin{rem} Note that the previous theorem cannot be formulated just using $\mathbb{O}_x(K)$ and $\mathbb{F}_x(K)$ (while Theorem~\ref{Knoroots} is formulated in this way). Indeed, if $L$ is first countable, then $\mathbb{F}_x(L)$ is singleton for each $x\in L$. Therefore $K$ and $K\times L$ cannot be distinguished using just the objects  $\mathbb{O}_x(K)$ and $\mathbb{F}_x(K)$ and there are first countable compact spaces of uncountable weight.
\end{rem}
\section{The ball of the Hilbert space, $P(A(\kappa))$ and $M(A(\kappa))$}

In the section we shall compute $\mathbb{O}_x$ for the ball of the
Hilbert space and its finite powers. In particular, we shall prove Theorem~\ref{OxBk} and Theorems~\ref{powerBall} and~\ref{productball} for the case $B(\kappa)$.
Recall that
$$B(\kappa) = \left\{(x_i)_{i<\kappa}\in \mathbb{R}^\kappa :
\sum_{i<\kappa}|x_i|^2\leq 1\right\}$$
endowed with the weak
topology of the Hilbert space $\ell_2(\kappa)$.
The weak topology clearly coincides with the pointwise one. We can identify this space by the
obvious homeomorphism with
$$B(\kappa) \approx \left\{(x_i)_{i<\kappa}\in \mathbb{R}^\kappa :
\sum_{i<\kappa}|x_i|\leq 1\right\}\subset\mathbb{R}^\kappa$$
with the pointwise topology. This compact is also
homeomorphic to the ball of $\ell_p(\kappa)$ for $1<p<\infty$ in
the weak topology and to the dual ball of $c_0(\kappa)$ in the
weak$^\ast$ topology. It is to be noticed that all the results
proved in this section hold true (with essentially identical
proof) if we substitute the space $B(\kappa)$ by $P(A(\kappa))
\approx \left\{(x_i)_{i<\kappa}\in [0,1]^\kappa :
\sum_{i<\kappa}x_i\leq 1\right\}$. The fiber orders of
$P(A(\kappa))$ will be computed again as one particular case of
our methods in spaces of probability measures. We shall also
notice that $P(A(\kappa))$ is not homeomorphic to the dual unit
ball of the Banach space of continuous functions $C(A(\kappa))$ in its weak$^\ast$ topology.

For a subset $M$ of $\kappa$, we consider $B(M) =
\left\{(x_i)_{i\in M}\in \mathbb{R}^M : \sum_{i\in M}|x_i|\leq
1\right\}$, and for $M\subset N$ we have the natural projection
$p_{MN}:B(M)\To B(N)$ given by $p((x_i)_{i\in M}) = (x_i)_{i\in N}$.
Thus every $B(M)$ can be seen as a quotient of $B(\kappa)$ through
the projection $p_{\kappa M}:B(\kappa)\To B(M)$, and all quotients
of this type for $M$ countably infinite subset of $\kappa$ constitute a
cofinal $\sigma$-semilattice of $\mathcal{Q}_\omega(B(\kappa))$,
as it easily follows from Lemma~\ref{spectralfactor}. Hence, the
$\sigma$-typical surjection of $B(\kappa)$ is of the form
$p_{MN}:B(M)\To B(N)$ and its fiber orders are computed in the following way:

\begin{lem}\label{ballfiberlemma}
Let $p_{MN}:B(M)\To B(N)$ be as above, $x\in B(N)$ and $y^1,y^2\in
p_{MN}^{-1}(x)$. Then $y^1\leq y^2$ if and only if $\sum_{i\in
M\setminus N}|y^1_i| \leq
\sum_{i\in M\setminus N}|y^2_i|$.
\end{lem}

Proof: Set $M^\ast =
M\setminus N$. Let $y$ be any point of $p_{MN}^{-1}(x)$. A basic
neighborhood of $y$ is of the form
$$U = \{z\in B(M) : z_i\in W_i\ \text{ for } i\in F\}$$
where $F$ is a finite subset of $M$ and $W_i$ is an open real
interval containing $y_i$, for every $i\in F$. Let $a_i =\inf\{|t|
: t\in W_i\}$ be the distance of the interval $W_i$ to 0. Then,
the image of the above typical basic neighborhood $U$ under
$p_{NM}$ is the following:

\begin{align*}p_{MN}(U) &= \left\{z\in B(N) : z_i\in W_i\ \text{ for } i\in F\cap N \right\},  \\ & \qquad\qquad\qquad\qquad\text{ if }0\in W_i \text{ for all } i\in F\cap M^\ast;\\
\ &\\
    p_{MN}(U) &= \smash{\left\{z\in B(N) : z_i\in W_i\ \text{ for } i\in F\cap N \phantom{\sum_{i\in M}}\right.}
    \\ & \qquad\qquad\qquad\qquad\text{ and }
\left.\sum_{i\in N}|z_i| < 1 - \sum_{i\in F\cap M^\ast}a_i\right\} \text{ otherwise.}
\end{align*}
This means that the images of the basic neighborhoods of $y$ are
the sets of the following form:
\begin{itemize}
\item If $y_i=0$ for all $i\in M^\ast$, then the images of the
basic neighborhoods of $y$ are the basic neighborhoods of $x$.

\item Otherwise, the images of basic neighborhoods of $y$ are the
sets of the form $V\cap \{z : \sum_{i\in N}|z_i|<1-r\}$ where $V$
is a basic neighborhood of $x$ and $r$ is any real number such
that $0\leq r <\sum_{i\in M^\ast}|y_i|$.
\end{itemize}

From this description, it is already clear that if $\sum_{i\in
M^\ast}|y^1_i| \leq \sum_{i\in M^\ast}|y^2_i|$ then $y^1\leq y^2$.
For the converse implication, it is enough to check that if
$r<s<1-\sum_{i\in N}|x_i|$ there is no neighborhood $V$ of $x$
such that $V\cap \{z : \sum_{i\in N}|z_i|<1-r\}\subseteq \{z :
\sum_{i\in N}|z_i|<1-s\}$. This follows from the fact that $N$ is
infinite: Suppose $V =\{z\in B(N) : z_i\in W_i \text{ for } i\in
F\}$ where $F$ is some finite subset of $N$ and $W_i$ are
intervals; take a number $1-s<t<1-r$ and $n\in N\setminus F$;
consider the element $y$ which agrees with $x$ on $F$, $y_n =
t-\sum_{i\in F}|x_i|$, and $y_i$ is 0 in all other coordinates. Then
$1-s<t=\sum_{i\in N}|y_i|<1-r$ and $y\in V$, so $y\in V\cap \{z :
\sum_{i\in
N}|z_i|<1-r\}\setminus \{z : \sum_{i\in N}|z_i|<1-s\}$.$\qed$

It follows from Lemma~\ref{ballfiberlemma} that
$\mathbb{O}_x(p_{MN})$ is order isomorphic to an interval $[a,b]$
if $\sum_{i\in N}|x_i|<1$, and $|\mathbb{O}_x(p_{MN})|=1$ if
$\sum_{i\in N}|x_i|=1$. From this, it is also clear that for $x\in
B(\kappa)$, $\mathbb{O}_x(B(\kappa))\cong [0,1]$ if
$\sum_{i<\kappa}|x_i|<1$, and $|\mathbb{O}_x(B(\kappa))|=1$ if
$\sum_{i<\kappa}|x_i|=1$. Theorem~\ref{OxBk} follows immediately now.

We notice that $B(\kappa)$ satisfies the hypotheses of both
Theorem~\ref{Knoroots} and Theorem~\ref{Knodecomposition}, taking
$[0,1]$ as irreducible ordered set. Every fiber of $p_{NM}$ has an
equivalence class which is a singleton, namely the class of the
minimum element, the one with $y_i=0$ fo all $i\in M^\ast$. This yields the proof of the case $B(\kappa)$ of Theorems~\ref{powerBall} and~\ref{productball}.

 We stated in the
introduction that whenever two points $x,y\in B(\kappa)^n$ have
the same number of coordinates in the sphere then there is a
homeomorphism $f:B(\kappa)^n\To B(\kappa)^n$ with $f(x)=y$. Let us
indicate why. It is enough to consider the case $n=1$. If
$\|x\|=\|y\|$ then we can find a linear isometry of the Hilbert
space $\ell_2(\kappa)$ onto itself sending $x$ to $y$. After this,
it remains to find some automorphism of $B(\kappa)$ sending some
element of norm $\lambda$ to some element of norm $\mu$ for every
$\lambda,\mu\in [0,1)$. View now again $B(\kappa) =
\{(x_i)_{i<\kappa} : \sum|x_i|\leq 1\}$. By the standard
homeomorphism, the norm function is transformed into $\|x\|^2 =
\sum|x_i|$. Consider an increasing homeomorphism $\phi:[-1,1]\To
[-1,1]$ such that $\phi(\lambda)=\mu$ and there exist the side
derivatives $\phi'_{+}(-1) = 1 = \phi'_{-}(1)$. Consider then
$f:B(\kappa)\To B(\kappa)$ given by $f((x_i)_{i<\kappa}) =
(y_i)_{i<\kappa}$ where $y_0=\phi(x_0)$ and $y_ i =
\frac{1-|\phi(x_0)|}{1-|x_0|}x_i$ for $i>0$. Notice that $f$ is a
homeomorphism and
$f(\lambda,0,0,\ldots) = (\mu,0,0,\ldots)$.

Let $M(K)$ denote the set of all Radon measures of variation at
most one (in other words, the dual ball of the Banach space
$C(K)$) endowed with the weak$^\ast$ topology. We know that
$\mathbb{O}_x(P(A(\kappa))$ is either a singleton or
order-isomorphic to $[0,1]$ for $x\in P(A(\kappa))$. A cofinal
$\sigma$-semilattice for $\mathcal{Q}_\omega(M(K))$ is formed by
the quotients of the form $M(p):M(K)\To M(L)$ where $p:K\To L$ is
an element of $\mathcal{Q}_\omega(K)$, and hence the
$\sigma$-typical surjection of $M(K)$ is of the form $M(g):M(X)\To
M(Y)$, where $g:X\To Y$ is the $\sigma$-typical surjection of $K$. Hence, in order to prove that $P(A(\kappa))\not\approx
M(A(\kappa))$, it is
enough to check the following:

\begin{prop}
Let $g:K\To L$ be continuous surjection and let $f=M(g):M(K)\To
M(L)$ be the induced map between the spaces of Radon measures of
variation at most 1. If there exists $x\in L$ such that $|\mathbb{O}_x(g)|>1$,
then $\mathbb{O}_0(f)$ is not linearly ordered.
\end{prop}

Proof: Take $y,z\in g^{-1}(x)$ such that $y\not\leq z$, so
that there is a neighborhood $U$ of $y$ such that $g(U)$
does not contain the image of any neighborhood of $z$, and moreover $z\not\in\overline{U}$. There is a net $(z_\alpha)$ in $K$ that converges to $z$
and with $g(z_\alpha)\not\in g(U)$ for every $\alpha$. Consider
the measures $\nu = \frac{1}{2}\delta_y - \frac{1}{2}\delta_{z}$
and $\mu=-\frac{1}{2}\delta_y+\frac{1}{2}\delta_z$. We claim that
these are two incomparable elements of $f^{-1}(0)$. We prove
that $\nu\not\leq\mu$ (that $\mu\not\leq\nu$ is done by analogy).
We consider $W=\{\lambda\in M(K) : \lambda(U)>\frac{3}{8}\}$,
which is a neighborhood of $\nu$. We claim that $f(W)$ does not
contain the $f$-image of any neighborhood of $\mu$. Notice that
$f(W)\subset \{\zeta\in M(L) : \zeta(g(U))>-\frac{1}{4}\}$.
Consider the measures $\mu_\alpha = -\frac{1}{2}\delta_y +
\frac{1}{2}\delta_{z_\alpha}$. Then $\mu_\alpha\To\mu$ and
$f(\mu_\alpha)(g(U)) = -\frac{1}{2}$ for each $\alpha$. In particular $f(\mu_\alpha)\notin f(W)$ for each $\alpha$. This witnesses that $\nu\not\le\mu$. $\qed$

\section{Computing images of neighborhoods in spaces of probability
measures}\label{sectLemma}

In order to compute the order of the fiber of a certain point
$y\in K_2$ in a surjection $f:K_1\To K_2$, we have to know how to
compute the images $f(U)$ of basic neighborhoods $U$ of points
$x\in f^{-1}(y)$. The surjections which appear in the cases that
we are going to study now are of the form $f=P(g):P(K)\To P(L)$
were $g:K\To L$ is a surjection between scattered compacta and $P$
is the functor of probability measure spaces. In this case, a
neighborhood basis of a measure $\mu\in P(K)$ is formed by the
sets of the following form:
$$\mathcal{U} = \{\nu \in P(K) : \nu(U_i)>c_i,\ i=1,\ldots,n\}$$
where $U_i$ are disjoint clopen subsets of $K$ from a given basis
of clopen sets, and $c_i$ are any numbers with $\mu(U_i)>c_i$. The
following lemma provides a computation of the image
$f(\mathcal{U})$ of such a neighborhood and will be applied
repeatedly in the future.

\begin{lem}\label{imageofneighborhood}
Let $g:K\To L$ be a surjection between compact spaces, let
$f=P(g):P(K)\To P(L)$, $U_1,\ldots,U_n$ be disjoint closed subsets
of $K$, $c_1,\ldots,c_n\geq 0$, and
$$\mathcal{U} = \{\nu\in P(K) : \nu(U_i)>c_i,\ i=1,\ldots,n\}.$$
Then
$$f(\mathcal{U}) = \left\{\lambda\in P(L) :
\lambda\left(g\left(\bigcup_{i\in A}U_i\right)\right)>\sum_{i\in
A} c_i\mbox{ for }A\subset \{1,\ldots,n\},
A\neq\emptyset\right\}.$$
\end{lem}

The fact that $f(\mathcal{U})$ is included in the righthand side
expression is trivial. The other inclusion is related to the
following numerical lemma:

\begin{lem}\label{numericalimage}
We consider numbers $c_1,\ldots,c_n\geq 0$, and $\alpha_A\geq 0$
for $A\subset\{1,\ldots,n\}$, $A\neq\emptyset$ such that for every
nonempty $A\subset\{1,\ldots,n\}$ we have that $\sum_{B\cap
A\neq\emptyset}\alpha_B > \sum_{i\in A}c_i$. Then for every
nonempty $A\subset\{1,\ldots,n\}$ there exist numbers
$\beta_{A,i}$, $i\in A$ such that $\sum_{i\in A}\beta_{A,i} =
\alpha_A$ and moreover, for every $i=1,\ldots,n$, $\sum_{A\ni
i}\beta_{A,i} > c_i$.
\end{lem}

Let us make some comment about the history of the lemmas. We first
had a long proof by induction of Lemma~\ref{numericalimage}. After
speaking about it with Richard Haydon, he indicated to us a more
elegant and shorter proof using combinatorial optimization that we
reproduce below. Later, after David Fremlin heard about it in the
Marczewski Centennial Conference in Bedlewo, he wrote a note
\cite{Fremlin} where he shows that actually
Lemma~\ref{imageofneighborhood} holds under more general
assumptions in $\mathcal K$-analytic spaces (our original
statement was only for scattered or metrizable compact sets).

We first notice how Lemma~\ref{imageofneighborhood} follows from
Lemma~\ref{numericalimage} in the cases when $K$ is either
scattered or metrizable, which is enough for the applications that
we present (for the general case we refer to \cite{Fremlin}).
Given a measure $\lambda$ in the righthand side of the conclusion
of Lemma~\ref{imageofneighborhood}, we consider $X_A = \cap_{i\in
A}g(U_i)\setminus\cup_{i\not\in A}g(U_i)$ and the numbers
$\alpha_A = \lambda(X_A)$ (note that each $X_A$ is Borel as it is
the difference of two closed sets), to which we can apply
Lemma~\ref{numericalimage} and obtain the numbers $\beta_{A,i}$.
We define a measure $\nu\in P(K)$ with $f(\nu)=\lambda$ in the
following way.

Suppose first that $K$ and $L$ are scattered, so that all Radon
measures on them are discrete, that is, determined by the measures
of singletons (in this case, we do not even need that the sets
$U_i$ are closed). If $y\in L\setminus\bigcup_1^n g(U_i)$ then we
pick a point $x_y\in g^{-1}(y)$ and declare $\nu(\{x_y\}) =
\lambda(\{y\})$. If $y\in X_A$ for some nonempty
$A\subset\{1,\ldots,n\}$ with $\alpha_A>0$, then we can choose elements $x_{y,i}\in
f^{-1}(y)\cap U_i$ for every $i\in A$ and then we declare
$\nu(\{x_{y,i}\}) = \frac{\beta_{A,i}}{\alpha_A}\lambda(\{y\})$.
In any other points $\nu(\{x\})=0$. This $\nu$ is a probability
measure on $K$ with $f(\nu)=\lambda$ and moreover
$\nu\in\mathcal{U}$ since $\nu(U_i)=\sum_{A\ni i}\beta_{A,i}$.

In the other case, suppose that $K$ and $L$ are metrizable. In
this case, Radon and Borel measures coincide.
 For each $A\subset\{1,\ldots,n\}$ and each $i\in A$,
we set $Y_{A,i}= g^{-1}(X_A)\cap U_i$ which is a nonempty Borel
set. We also define $X_\emptyset = L\setminus\bigcup_1^n g(U_i)$
and $Y_\emptyset = g^{-1}(X_\emptyset)$. By the
Jankov-Von~Neumann Uniformization Theorem \cite[Theorem 18.1]{Kechris},
there exists a measurable selection $s_{A,i}:X_A\To Y_{A,i}$ for
the the inverse of $g_{|Y_{A_i}}$ and also a measurable selection $s_\emptyset:X_\emptyset\To Y_\emptyset$ of the inverse of $g_{|Y_\emptyset}$. Consider
$$\nu = s_\emptyset(\lambda|_{X_\emptyset}) + \sum_{i\in
A}\frac{\beta_{A,i}}{\alpha_A}\cdot s_{A,i}(\lambda|_{X_{A,i}}).$$
Then $f(\nu)=\lambda$, and $\nu(U_i)=\sum_{A\ni i}\beta_{A,i}$, so
$\nu\in\mathcal{U}$.

For the proof of Lemma~\ref{numericalimage}, we shall use the so
called max-flow min-cut theorem, Theorem~\ref{minmax} below, from
combinatorial optimization. This result is originally due to Ford
and Fulkerson~\cite{maxflowmincut1} and Dantzig and
Fulkerson~\cite{maxflowmincut2}, and can be found in the
book~\cite[Theorem 10.3]{combinatorialoptimization}. We have to
recall some concepts from this area. A directed graph
(\emph{digraph} for short) is a couple $G=(V,A)$ where $V$ is a
finite set whose elements are called \emph{vertices}, and
$A\subset V\times V$ is a set whose elements are called \emph{arcs}.
An $s-t$-flow is a function $f:A\To(0,+\infty)$ which satifies the flow conservation law at
all points except $s$ and $t$:

$$\sum_{(x,u)\in A}f(x,u) = \sum_{(u,x)\in A}f(u,x) \text{\qquad for every }u\in V\setminus\{s,t\}.$$

In words, the flow entering $u$ equals the flow leaving $u$. The
value of the flow $f$ is the net amount of flow leaving $s$, which
happens to be equal to the net amount of flow entering $t$,

$$\mbox{value}(f) = \sum_{(s,x)\in A}f(s,x) - \sum_{(x,s)\in A}f(x,s) = \sum_{(x,t)\in A}f(x,t) -\sum_{(t,x)\in A}f(t,x).$$

Let us a consider a function $c:A\To(0,+\infty)$ that we call a
capacity function. A flow $f$ is said to be under $c$ if
$f(u,v)\leq c(u,v)$ for every $(u,v)\in A$. Given a set $B\subset
A$, the capacity of $B$ is $c(B) = \sum_{(u,v)\in B}c(u,v)$.

For a subset $U\subset V$, we denote by $\delta(U)$ the set of all
arcs which leave $U$ and enter $V\setminus U$, that is,
$$\delta(U) = \{(u,v)\in A : u\in U, v\not\in U\}.$$

For $s,t\in V$, an $s-t$ cut is a set of arcs of the form
$\delta(U)$, where $U\subset V$ with $s\in U$ and $t\not\in U$.

\begin{thm}[max-flow min-cut theorem]\label{minmax} Let $G=(V,A)$
be a digraph, $t,s\in V$ and $c:A\To\mathbb{R}_+$ a capacity
function. Then the maximum value of an $s-t$ flow under $c$ equals
the minimum capacity of an $s-t$-cut,
$$\max\{\mbox{\rm value}(f) : f\leq c \text{ is an }s-t\text{ flow}\} =
\min\{c(\delta(U)) : U\subset V, s\in U, t\not\in U\}.$$
\end{thm}

Proof of Lemma~\ref{numericalimage}: We shall denote by
$\mathcal{P}_n$ the family of all nonempty subsets of
$\{1,\ldots,n\}$. First we consider numbers $c_i'>c_i$ for every
$i\leq n$ such that the inequalities
$$\sum_{B\cap A\neq\emptyset}\alpha_B > \sum_{i\in A}c'_i$$ still
hold. We consider a digraph $G=(V,A)$ where the set of vertices is
$$V = \{s\}\cup \{p_i : 1\leq i\leq n\}\cup \{q_A : A\in
\mathcal{P}_n\}\cup \{t\},$$ and the set of arcs is
$$A = \{(s,p_i) : 1\leq i\leq n\}\cup \{(p_i,q_A) : i\in A\}\cup
\{(q_A,t) : A\in\mathcal{P}_n\}.$$ Let $M\in(0,+\infty)$ be such
that $M>\sum_{i=1}^n c'_i$. We define a capacity function
$c:A\To\mathbb{R}_+$ as
\begin{itemize}
\item $c(s,p_i) = c'_i$. \item $c(p_i,q_A)=M$. \item $c(p_i,q_A) =
\alpha_A$.
\end{itemize}

\emph{Claim A}: The minimal capacity of an $s-t$ cut in $G$ equals
$\sum_{i=1}^n c'_i$.

\smallskip

Proof of Claim A: If $U=\{s\}$, then $c(\delta(U)) =
\sum_1^n c'_i$. We suppose that $\delta(U)$ is an arbitrary $s-t$
cut and we show that its capacity is larger than $\sum_1^n c'_i$.
Let $A = \{i\leq n : p_i\in U\}$. If there exists $B\in
\mathcal{P}_n$ such that $q_B\not\in U$ and $A\cap
B\neq\emptyset$, then there exists $(p_i,q_B)\in \delta(U)$ and in
particular $c(\delta(U))\geq c(p_i,q_B)=M>\sum_1^n c_i'$.
Hence, we can suppose that $q_B\in U$ whenever $A\cap
B\neq\emptyset$, therefore $(q_B,t)\in \delta(U)$ whenever $A\cap
B\neq\emptyset$, and $$c(\delta(U)) \geq \sum_{i\not\in A}c(s,p_i)
+ \sum_{B\cap A\neq\emptyset}c(q_B,t) = \sum_{i\not\in A}c'_i +
\sum_{B\cap A\neq\emptyset}\alpha_B \geq \sum_{i\not\in A}c'_i +
\sum_{i\in A}c'_i = \sum_{i=1}^n c'_i.$$
\

By Claim $A$ and Theorem~\ref{minmax}, there exists an $r-s$ flow
$f\leq c$ and value equal to $\sum_1^n c'_i$. Notice that
$f(s,p_i)\leq c'_i$ but $\sum_1^n f(s,p_i) = \mbox{value}(f) = \sum_1^n
c'_i$, hence $f(s,p_i) = c'_i$.

By the flow conservation law at the vertex $q_A$, for every
$A\in\mathcal{P}_n$ we have that $$\sum_{i\in A}f(p_i,q_A) =
f(q_A,t)\leq c(q_A,t)= \alpha_A,$$ therefore we can choose numbers
$\beta_{A,i}$ for $i\in A$ such that $\beta_{A,i}\geq f(p_i,q_A)$
and $\sum_{i\in A}\beta_{A,i}=\alpha_A$. We claim that these
numbers have the desired property. To check this, we use again the
flow conservation law now at a vertex $p_i$, \[c_i<c'_i = f(s,p_i)
= \sum_{A\ni i}f(p_i,q_A) \leq \sum_{A\ni i}\beta_{A,i}.\qed \]

\section{Fiber orders of the probability measures on a scattered compact}

As we already mentioned, it is a standard fact that if $K$ is a
totally disconnected compact and $\mathcal{B}$ is a basis for the
topology of $K$ consisting of clopen sets, then a basis for the
topology of $P(K)$ consists of the sets of the form $\{\mu\in P(K)
: \mu(U_i)>c_i : i=1,\ldots,n\}$, where the $c_i$'s are positive
numbers and the $U_i$'s are disjoint basic clopen sets. When $K$
is scattered, all measures from $P(K)$ are discrete, and this
allows to find a finer neighborhood basis which will be quite
useful for
us. To avoid heavy notation, we write $\mu(x)$ instead of $\mu(\{x\})$ to denote the measure of a singleton.

\begin{lem}
Let $K$ be a scattered compact space, $\mu\in P(K)$ and let
$\mathcal{B}$ be a basis of the topology of $K$ consisting of
clopen sets. A neighborhood basis of $\mu$ consists of the sets of
the form $\{\nu:\nu(U_i)>c_i\mbox{ for }i=1,\dots,n\}$, where
$U_1,\dots,U_n$ are pairwise disjoint basic clopen neighborhoods of
points $x_1$,$\ldots$,$x_n$ of $K$ and $c_1,\dots,c_n$ are
positive numbers with $\mu(x_i)>c_i$.
\end{lem}

Proof: Consider a neigborhood of $\mu$ of the form $V=\{\nu\in
P(K) : \nu(V_j)>d_j, j=1,\ldots,n\}$ for some disjoint basic
clopen neighborhoods $V_j$ with $\mu(V_j)>d_j$. Since $\mu$ is
discrete, for every $j$ we can find a finite family of points
$\{x^j_i : i\in F_j\}$ such that $\sum\{\mu(x^j_i) : i\in
F_j\}>d_j$. Consider numbers $d_i^j<\mu(x^j_i)$ such that
$\sum_{i\in F_j}d_i^j>d_j$, and disjoint basic clopen sets
$x^j_i\in U^j_i\subset V_j$, then \[\mu\in\{\nu\in P(K) :
\nu(U^j_i)>d_i^j, j=1,\ldots,n,i\in F_j\}\subset V \qed\]

For the rest of the section, we fix $g:K\To L$ to be a surjection
between scattered compact spaces and $f=P(g):P(K)\To P(L)$ the
induced map between the spaces of probability measures. Note that
the norm we use below is the $\ell_1$-norm, i.e.

$$\|\mu-\nu\|=\sum_{k\in K} |\mu(k)-\nu(k)|\mbox{ for }\mu,\nu\in P(K).$$

\begin{thm}\label{sumofprobabilities}
Let $\mu = \sum_{i\in I} r_i\delta_{x_i}$ be a probability measure on $L$, where $I=\nat$ or $I=\{1,\dots,N\}$ for some $N\in\nat$, $x_i$, $i\in I$ are pairwise distinct points in $L$ and $r_i>0$ for all $i\in I$.
Then, the natural bijection $\prod_{i\in I}f^{-1}(\delta_{x_i})\To
f^{-1}(\mu)$ given by $(\nu_i)_{i\in I} \mapsto \sum_{i\in I} r_i\nu_i$ is an
order-isomorphism. In particular $\mathbb{F}_\mu(f)\cong \prod_{i\in I}
\mathbb{F}_{\delta_{x_i}}(f)$ and $\mathbb{O}_\mu(f)\cong \prod_{i\in I}
\mathbb{O}_{\delta_{x_i}}(f)$.
\end{thm}

Proof: Consider the mapping $\Phi: P(K)^I\to P(K)$ defined
by
\begin{equation*}
\Phi((\nu_i)_{i\in I})=\sum_{i\in I} r_i\nu_i.
\end{equation*}
It is easy to check that $\Phi$ is a continuous surjection. Moreover, as it is affine,
it maps $\prod_{i\in I}f^{-1}(\delta_{x_i})$ bijectively onto $f^{-1}(\mu)$. We will show that
the restriction of $\Phi$ to $\prod_{i\in I}f^{-1}(\delta_{x_i})$ is an order-isomorphism.

\smallskip

First, suppose that $\sum_{i\in I} r_i\nu_i \leq \sum_{i\in I} r_i \nu'_i$
with $\nu_i,\nu'_i\in f^{-1}(\delta_{x_i})$ for $i\in I$
and we shall prove that $\nu_1\leq \nu'_1$. We consider a typical neighborhood of $\nu_1$
of the form
$$U = \{\nu\in P(K):\nu(U_j)>c_j\mbox{ for
}j=1,\dots,n\}$$
where each $U_j$ is a clopen neighborhood of some
$a_j$ satisfying $\nu_1(a_j)>c_j$ and $U_j$'s are pairwise disjoint.
We want to find a neighborhood $V$ of
$\nu'_1$ such that $f(V)\subset f(U)$.

We consider a very small number $\varepsilon>0$, namely such that
\begin{equation*}
\varepsilon<r_1\mbox{\quad and\quad} (r_1+\varepsilon)(c_j+2\varepsilon)<r_1\nu_1(a_j)
\mbox{ for }j=1,\dots,n.
\end{equation*}
Choose $k\in I$ such that
$$\sum_{i>k}r_i<\frac{1}{4}r_1\varepsilon(r_1-\varepsilon)$$
and $H$ a clopen subset of $L$ such that $x_1\in H$ but $x_i\not\in H$ for
$2\leq i\leq k$.
The following is a neighborhood of $\sum_{i\in I}
r_i\nu_i$:
$$U^0 = \{\nu\in P(K):\nu(U_j\cap g^{-1}(H))>(r_1+\varepsilon)(c_j+2\varepsilon) \mbox{ for }j=1,\dots,n\}$$
By our assumption, there exists $V^0$ a neighborhood of $\sum_{i\in I}
r_i\nu'_i$ such that $f(V^0)\subset f(U^0)$. We take $V^0$ to be
of the typical form
$$V^0 = \{\nu\in P(K):\nu(V_{j})>d_{j}\mbox{ for
}j\in J\}$$
where $J$ is a finite set, $V_j$ is a clopen
neighborhood of some $b_j$ satisfying
$$\sum_{i\in I} r_i\nu'_i(b_j)>d_j$$
and $V_j$'s are pairwise disjoint. We let
$$J_i = \{j\in J : g(b_j)=x_i\}.$$
Without loss of generality we suppose
that $V_j\subset g^{-1}(H)$ for $j\in J_1$, and $V_j\cap
g^{-1}(H)=\emptyset$ for $j\in J_2\cup\cdots\cup J_k$. Notice that
\begin{equation}\label{remainderdij}
\sum\left\{d_j : j\in
\bigcup_{i>k}J_i\right\}<\frac{1}{4}r_1\varepsilon(r_1-\varepsilon).
\end{equation}
Consider now
$$V = \{\nu\in P(K):\nu(V_{j})>d_{j}/r_1\mbox{ for }j\in J_1\}\cap\{\nu\in P(K) : \nu(K\setminus g^{-1}(H))<\varepsilon/2\}.$$
This is a neighborhood of $\nu'_1$ (notice that $\nu'_1(K\setminus g^{-1}(H))=0$
and $\nu'_1(V_j)\ge\nu'_1(b_j)>d_j/r_1$ for $j\in J_1$).
We claim that $f(V)\subset f(U)$.

So take $\xi_1\in V$. We can easily find $\xi_2\in V$ with
$\|\xi_2-\xi_1\|<\varepsilon$ such that $\xi_2(K\setminus
g^{-1}(H))=0$. We pick a measure $\lambda\in P(K)$ with
$\lambda(g^{-1}(H))=0$ and $\lambda(V_j)>d_j/(1-r_1)$ for $j\in
J_2\cup\cdots\cup J_k$. Then the measure $\xi_3 = r_1\xi_2 +
(1-r_1)\lambda$ satisfies $\xi_3(V_j)>d_j$ for $j\in
J_1\cup\cdots\cup J_k$. By (\ref{remainderdij}) we may find
$\xi_4\in V^0$ such that
$$\|\xi_4-\xi_3\|<\frac{1}{2}r_1\varepsilon(r_1-\varepsilon).$$
Set
$$r=\xi_4(g^{-1}(H)) \mbox{ and } \xi_5=\frac1r{\xi_4}_{|g^{-1}(H)}.$$
We have
\begin{equation}\label{r1-r}
|r_1-r| = | \xi_3(g^{-1}(H))-
\xi_4(g^{-1}(H))|<\frac{1}{2}r_1\varepsilon(r_1-\varepsilon)<\varepsilon
\end{equation}
and
\begin{align*}
\|\xi_2-\xi_5\| &= \sum_{t\in g^{-1}(H)}|\xi_2(t)-\xi_5(t)| =
\sum_{t\in g^{-1}(H)}\left|\frac{\xi_3(t)}{r_1}-\frac{\xi_4(t)}{r}\right|\\
&= \sum_{t\in g^{-1}(H)}\frac{1}{r_1 r}|r\xi_3(t)-r_1\xi_4(t)| \leq \frac{1}{r_1
r}\|r\xi_3-r_1\xi_4\|\\
&\leq  \frac{1}{r_1
r}\|(r_1-r)\xi_4\| + \frac{1}{r_1 r}\|r\xi_3-r\xi_4\|\\
&\le\frac1{r_1 r}|r_1-r|+\frac1{r_1}\|\xi_3-\xi_4\| \\
&< \frac{\varepsilon(r_1-\varepsilon)}{2 r} +
\frac{1}{2}\varepsilon(r_1-\varepsilon)<\frac{\varepsilon}{2}+\frac{\varepsilon}{2}
= \varepsilon.
\end{align*}
The first inequality on the last line follows from the first inequality of (\ref{r1-r}),
for the second one we use the fact that $r_1-\varepsilon<r$.
It follows that $\|\xi_1-\xi_5\|<2\varepsilon$, and hence
\begin{equation}\label{perturbation}
\|f(\xi_1)-f(\xi_5)\|<2\varepsilon
\end{equation}
as well.

Now, $\xi_4\in V_0$, hence $f(\xi_4)\in f(V_0)\subset f(U_0)$.
By the description of $f(U^0)$ given by
Lemma~\ref{imageofneighborhood}, the fact that all
clopen subsets of $K$ appearing in the definition of $U^0$ are
contained in $g^{-1}(H)$ implies that
\begin{equation}\label{perturbed}
\begin{aligned} f(\xi_5) &\in f(\{ \nu\in
P(K):\nu(U_j\cap g^{-1}(H))>(r_1+\varepsilon)(c_j+2\varepsilon)/r
\mbox{ for }j=1,\dots,n \})\\ &\subset f(\{ \nu\in P(K):\nu(U_j)>c_j+2\varepsilon \mbox{ for }j=1,\dots,n \}).
\end{aligned}
\end{equation}
The inclusion above follows from (\ref{r1-r}) -- note that $r_1+\varepsilon>r$.
Finally, using (\ref{perturbed}) and (\ref{perturbation}) it easily follows from
 Lemma~\ref{imageofneighborhood} that $f(\xi_1)\in f(U)$ which completes the proof of the first implication.

\smallskip

We pass now to the converse implication. So we assume that
$\nu_i\leq\nu'_i$ for every $i$, and we shall see that $\sum_{i\in I}
r_i\nu_i\leq \sum_{i\in I} r_i\nu'_i$.

Let $U$ be a neighborhood of $\sum r_i\nu_i$ in $P(K)$.
By the continuity of $\Phi$ and the definition of the product topology there is some $k\in I$
and neigborhoods $U_i$ of $\nu_i$ for $i\le k$ such that
\begin{equation}\label{inclusion}
\left\{ \sum_{i\in I}r_i\lambda_i : \lambda_i\in P(K)\mbox{ for }i\in I,
\lambda_i\in U_i\mbox{ for }i\le k\right\}\subset U.
\end{equation}
As $\nu_i\le\nu'_i$ for all $i\in I$, there is, for each $i\le k$, a neighborhood $V_i$ of $\nu'_i$ such that $f(V_i)\subset f(U_i)$.

Now we are going to specify the form of $V_i$'s. Let $H_1,\dots,H_k$ be pairwise disjoint clopen subsets of $L$ containing $x_1,\dots,x_k$, respectively. Then we can without loss of generality suppose that for each $i\le k$ we have
$$V_i=\{\lambda\in P(K): \lambda(V_i^j)>d_i^j \mbox{ for }j\in J_i\}$$
where $J_i$ is a finite set, $d_i^j>0$ for $j\in J_i$ and $V_i^j$, $j\in J_i$, are pairwise disjoint clopen subsets of $g^{-1}(H_i)$. Set
$$V=\{\lambda\in P(K): \lambda(V_i^j)>r_i d_i^j \mbox{ for }j\in J_i, i\le k\}.$$
Then $V$ is clearly a neighborhood of $\sum_{i\in I} r_i\nu'_i$. We claim that $f(V)\subset f(U)$.

Let $\lambda\in V$ be arbitrary. Choose $\delta>0$ such that
\begin{alignat}{3}
(1+\delta)\sum_{j\in J_i} d_i^j&<1, & i&=1,\dots,k; \label{suma}\\
\lambda(V_i^j)&>(1+\delta)r_i d_i^j, &\quad j&\in J_i,i=1,\dots,k.\label{Vij}
\end{alignat}
Further, define the following measures:
\begin{align*}
\sigma_i&=\sum_{j\in J_i} (1+\delta)d_i^j\frac{\lambda_{|V_i^j}}{\lambda(V_i^j)}, \qquad i=1,\dots,k, \\
\tau&=\lambda-\sum_{i=1}^k r_i\sigma_i.
\end{align*}
All $\sigma_i$'s are clearly positive measures. Moreover, $\tau$ is positive, too, by (\ref{Vij}) as
$$\tau_{|V_i^j}=\lambda_{|V_i^j}\left(1-\frac{(1+\delta) r_i d_i^j}{\lambda(V_i^j)}\right)$$
for $j\in J_i$, $i=1,\dots,k$.

It follows from (\ref{suma}) that $\sigma_i(K)<1$ for each $i=1,\dots,k$, and so $\tau(K)>0$.
For $i=1,\dots,k$ set
$$\theta_i=\sigma_i+\frac{1-\sigma_i(K)}{\tau(K)}\tau.$$
Then $\theta_i\in P(K)$. Moreover,
$$\theta_i(V_i^j)\ge\sigma_i(V_i^j)=(1+\delta)d_i^j>d_i^j$$
for $j\in J_i$, hence $\theta_i\in V_i$ for $i=1,\dots,k$. We
claim that
\begin{equation}\label{combination}
\lambda\in \left\{\sum_{i\in I} r_i\lambda_i: \lambda_i\in V_i\mbox{ for }i=1,\dots,k\right\}.
\end{equation}
Indeed, we can take $\lambda_i=\theta_i$ for $i=1,\dots,k$. To see
this, we have to check that
$$\vartheta=\lambda-\sum_{i=1}^k r_i\theta_i$$
is a nonnegative measure. Namely,
\begin{align*}
\vartheta &= \lambda-\sum_{i=1}^k r_i\theta_i = \lambda - \sum_{i=1}^k r_i\left(\sigma_i + \frac{1-\sigma_i(K)}{\tau(K)}\tau\right)\\
&= \lambda - \sum_{i=1}^k r_i \sigma_i - \sum_{i=1}^k r_i
\frac{1-\sigma_i(K)}{\tau(K)}\tau = \tau - \sum_{i=1}^k r_i
\frac{1-\sigma_i(K)}{\tau(K)}\tau\\
&= \frac{\tau}{\tau(K)}\left(\tau(K)-\sum_{i=1}^k r_i +
\sum_{i=1}^k r_i\sigma_i(K)\right),
\end{align*}
which is positive because $\sum_1^k r_i \le 1 = \lambda(K) = \tau(K)
+ \sum_1^k r_i\sigma_i(K)$. Thus, (\ref{combination}) is proved.
Using (\ref{combination}) and (\ref{inclusion}) we get by Lemma~\ref{imageofneighborhood} that $f(\lambda)\in f(U)$ which completes the proof.$\qed$

Let $g:K\To L$ be a continuous surjection, $x\in L$ and
$y_1$,$\ldots$,$y_n$ elements of the fiber $g^{-1}(x)$, we define
$\langle y_1,\ldots, y_n\rangle$ to be the set of all elements
$z\in g^{-1}(x)$ such that for every neighborhoods $U_1$,$\ldots$,
$U_n$ of $y_1$,$\ldots$,$y_n$ respectively there exists a
neighborhood $V$ of $z$ such that $g(V)\subset
g(U_1)\cup\cdots\cup g(U_n)$.

Notice some elementary properties, for instance $\langle y\rangle
= \{z\in g^{-1}(x) : z\geq y\}$ and $\langle Y\rangle\subset
\langle Y'\rangle$ whenever $Y\subset Y'$. The
$\langle\cdot\rangle$-operation provides in general a finer
structure on the fiber $g^{-1}(x)$ than the one given by the
order, and it is needed to determine the fiber order on spaces
$P(K)$ in terms of
the fibers of $K$. To avoid heavy notation, for a measure $\nu$, we often write $\nu\langle\cdot\rangle$ and $\nu\{\cdot\}$ instead of $\nu(\langle\cdot\rangle)$ and $\nu(\{\cdot\})$.

\begin{thm}\label{orderDirac}
Let $\nu$,$\nu'$ be elements of $f^{-1}(\delta_x)$. Then
$\nu\leq\nu'$ if and only if for every elements
$y_1$,$\ldots$,$y_n$ of $g^{-1}(x)$, $\nu\langle y_1,\ldots, y_n\rangle\leq\nu'\langle y_1,\ldots, y_n\rangle$.
\end{thm}

Proof: Suppose first that $\nu\leq\nu'$, and let
$y_1$,$\ldots$,$y_n$ be elements of $g^{-1}(x)$. If $\nu\langle
y_1,\ldots, y_n\rangle> \nu'\langle y_1,\ldots, y_n\rangle$, this
would mean that we can find elements $u_1$,$\ldots$,$u_r$ in
$\langle y_1,\ldots, y_n\rangle$, and elements
$v_1$,$\ldots$,$v_s$ in $g^{-1}(x)\setminus\langle y_1,\ldots,
y_n\rangle$, and a number $\xi>0$ such that $\sum_1^r\nu(u_i)>\xi$
and $\sum_1^s\nu'(v_i)>1-\xi$. Since $v_j\not\in\langle
u_1\ldots,u_r\rangle$ for any $j$ we can find neighborhoods
$U_{ij}$ of $u_i$ such that $\bigcup_{i=1}^r g(U_{ij})$ does not
contain the $g$-image of any neighborhood of $v_j$. Call
$U_i=\cap_{j=1}^s U_{ij}$. Consider a neighborhood of $\nu$ of the
form
$$U = \{\lambda\in P(K) : \lambda(U_i)> d_i,
i=1,\ldots,r\}$$ where $d_i<\nu(u_i)$ and $\sum_{1}^r d_i>\xi$. We
claim that $f(U)$ does not contain the image of any neighborhood
of $\nu'$, contradicting the fact that $\nu\leq\nu'$. So take any
neighborhood of $\nu'$, that we can suppose of the normal form
$$V = \{\lambda\in P(K) : \lambda(V_j)>e_j, j=1,\ldots,k\}$$
where the $V_j$'s are disjoint neighborhoods of points $w_j$ with
$\nu'(w_j)>e_j$ and moreover we can assume that $w_j=v_j$ for
$j=1,\ldots,s$. For every $j=1,\ldots,s$ we can find a point
$x_j\in g(V_j)\setminus\bigcup_{i=1}^r g(U_i)$. Consider the
measure $\lambda\in P(L)$ such that
$$\lambda(x_j) = \sum\{ \nu'(v_{j'}) : x_{j'}=x_j\}\text{ for }
j=1,\ldots,s$$
and $\lambda(x) = 1-\sum_{j=1}^s\nu'(v_j)$. Using
Lemma~\ref{imageofneighborhood} it easily follows that $\lambda\in
f(V)$. Further notice that
$$\lambda\left(\bigcup_{i=1}^r g(U_i)\right)\leq 1- \sum_{i=1}^s\lambda(x_i) =
1-\sum_{i=1}^s\nu'(v_i) < \xi.$$
This implies that $\lambda\not\in f(U)$ because otherwise we
should have that
$$\lambda\left(\bigcup_{i=1}^r g(U_i)\right)>\sum_1^r d_i
>\xi.$$

 Now we suppose that $\nu \langle Y\rangle \leq \nu'\langle
Y\rangle$ for every finite set $Y\subset g^{-1}(x)$. We want to
see that $\nu\leq\nu'$ so we take a typical neighborhood of $\nu$
of the form
$$U=\{\lambda\in P(K) : \lambda(U_i)>c_i, i=1,\ldots,n\}$$
where
the $U_i$'s are disjoint clopen neighborhoods of points $y_i$ such
that $\nu(y_i)>c_i$. For every nonempty $A\subset\{1,\ldots,n\}$
we have
$$\nu'\langle y_i : i\in A\rangle \geq \nu\langle y_i : i\in
A\rangle\geq \nu\{y_i:i\in A\}>\sum_{i\in A}c_i,$$
and so there exists a finite set of points $\{z_j : j\in
F_A\}\subset\langle y_i : i\in A\rangle$ such that
$$\sum_{j\in F_A} \nu'(z_j)>\sum_{i\in A}c_i.$$
Pick numbers $\xi_j<\nu'(z_j)$
such that $\sum_{j\in F_A}\xi_j>\sum_{i\in A}c_i$. For every $j\in
F_A$, since $z_j\in\langle y_i : i\in A\rangle$ we can find a
clopen neighborhood $V_j$ of $z_j$ such that $g(V_j)\subset \bigcup_{i\in
A}g(U_i)$. We can suppose that $V_j\cap V_{j'}=\emptyset$ for
different $j,j'\in F_A$. Now, for every $A$ the following is a
neighborhood of $\nu'$:
$$V^A = \{\lambda\in P(K) : \lambda(V_j)>\xi_j, j\in F_A\}$$
Let $V = \bigcap\{V^A : \emptyset\neq A\subset\{1,\ldots,n\}\}$.
We claim that $f(V)\subset f(U)$. Take $\lambda\in f(V)$.
According to Lemma~\ref{imageofneighborhood} we have to check that
for every nonempty $A\subset\{1,\ldots,n\}$,
$\lambda(\bigcup_{i\in A} g(U_i))>\sum_{i\in A}c_i$. Since
$\lambda\in f(V)\subset f(V^A)$, by the same lemma we know that
$\lambda(\bigcup_{j\in F_A} g(V_j))>\sum_{j\in
F_A}\xi_j>\sum_{i\in A}c_i$, and on the other hand $\bigcup_{j\in
F_A}g(V_j)\subset
\bigcup_{i\in A}g(U_i)$.$\qed$

We finish this section by the following proposition which will enable us to verify the assumptions of Theorems~\ref{Knoroots} and~\ref{Knodecomposition} for spaces of the form $P(K)$.

\begin{prop}\label{assumdecomp}
Let $K$ be a scattered compact space and $O$ a connected irreducible ordered set with $|O|>1$.
\begin{itemize}
    \item[(i)] Suppose that there is a point $x\in K$ such that $\mathbb{O}_{\delta_x}(P(K))\cong O$ and for each $y\in K\setminus \{x\}$ we have $\mathbb{O}_{\delta_y}(P(K))\not\cong O^k$ for any $k\ge 1$. Then $P(K)$ satisfies the assumption of Theorem~\ref{Knoroots}.
    \item[(ii)] Suppose that for the $\sigma$-typical surjection $f:L\to M$ (where $L$ and $M$ are metrizable quotients of $K$) there is $x\in M$ such that $\mathbb{O}_{\delta_x}(P(f))\cong O$, $\mathbb{F}_{\delta_x}(P(f))$ has one equivalence class which is a singleton and $\mathbb{F}_{\delta_y}(P(f))$ is a singleton for each $y\in M\setminus\{x\}$. Then $P(K)$ satisfies the assumptions of Theorem~\ref{Knodecomposition}.
\end{itemize}
\end{prop}

Proof: (i) We have $\mathbb{O}_{\delta_x}(P(K))\cong O$. Further,
suppose that $\mathbb{O}_\mu(P(K))\cong O^k$ for some $k>1$ for some $\mu\in P(K)$.
Let $C$ be a countable set supporting $\mu$. Then
it follows from Theorem~\ref{sumofprobabilities} and
Corollary~\ref{productofirreducibles} that for the $\sigma$-typical surjection $f$ of $K$ there is some $y\in C\setminus\{x\}$ such that $\mathbb{O}_{\delta_y}(P(f))\cong O^j$ for some
$j\ge 1$. Now, as $C$ is countable, it implies that there is $y\in C^\ast = C\setminus\{x\}$ such that in each cofinal $\sigma$-semilattice in $K$ there is a surjection $f$ such that $\mathbb{O}_{\delta_y}(P(f))\cong O^j$ for some $j\ge1$, which contradicts our assumptions. (Otherwise, for every $y\in C^\ast$ there would be a cofinal $\sigma$-lattice $\mathcal{S}_y\subset\mathcal{Q}_\omega(K)$ with $\mathbb{O}_{\delta_y}(P(f))\not\cong O^j$, for every $j$; an obvious improvement of Theorem~\ref{ourspectraltheorem} shows that $\bigcap_{y\in C^\ast}\mathcal{S}_k$ is a cofinal $\sigma$-semilattice which leads to a contradiction). Thus we have verified
the assumptions of Theorem~\ref{Knoroots}.

\smallskip

(ii) Consider $f:L\to M$ and $x\in M$ as in the assumptions. Clearly
$\mathbb{F}_{\delta_y}(P(f))$ is a singleton for each $y\in M\setminus\{x\}$.
Therefore we get, by
Theorem~\ref{sumofprobabilities} that $\mathbb{O}_\mu(P(f))\cong
O$ if $\mu(x)>0$ and $\mathbb{F}_\mu(P(f))$ is a singleton if
$\mu(x)=0$. In this way we have verified conditions (1)--(3) of
Theorem~\ref{Knodecomposition}. The remaining condition (4)
follows immediately from Theorem~\ref{sumofprobabilities}.$\qed$

\section{Examples of spaces of probability measures}

\subsection{The space $\sigma_n(\kappa)$} In this section we are going to prove Theorem~\ref{symmetric}, the case of $P(\sigma_n(\kappa))$ of Theorem~\ref{powerBall} and the case of $P(A(\kappa))$ of Theorem~\ref{productball}.

For $N\subset
M$, let $g_{MN}:\sigma_n(M)\To \sigma_n(N)$ be the continuous
surjection given by $g=g_{NM}(x) = x\cap N$. The $\sigma$-typical
surjection of $P(\sigma_n(\kappa))$ is of the form $f=P(g):
P(\sigma_n(M))\To P(\sigma_n(N))$ for $M\subset N$ infinite
countable subsets of $\kappa$ such that $M^\ast = M\setminus N$ is
infinite. The computation of the fiber order and the
$\langle\cdot\rangle$-operation is done as follows:

For $x\in \sigma_n(N)$, $g^{-1}(x) = \{x\cup y : y\subset M^\ast,
|y|\leq n-|x|\}$. A basic neighborhood of such $x\cup y\in
g^{-1}(x)$ is of the form $$U = \{z\in\sigma_n(M) : x\cup y
\subset z, z\cap u=\emptyset, z\cap v=\emptyset\},$$ where
$u\subset M^\ast$ and $v\subset N$ are finite sets. The image of
such a neighborhood equals
$$g(U) = \{z\in\sigma_n(N): x\subset z, |z|\leq n-|x\cup y|, z\cap v=\emptyset\}.$$
From this it is clear that for $w,w'\in g^{-1}(x)$ we have that
$w\leq w'$ if and only if $|w|\leq |w'|$ and also that
$$\langle
w_1,\ldots,w_k\rangle = \{ w\in g^{-1}(x) : |w|\geq
\min(|w_1|,\ldots,|w_k|)\}.$$
Thus, if we go now to the spaces of probabilities, for
$f=P(g):P(\sigma_n(M))\To P(\sigma_n(N))$, for two measures
$\nu,\nu'\in f^{-1}(\delta_x)$ we have that $\nu\leq \nu'$ if and
only if for every $k=1,\ldots,n-|x|$ we have that
$$\nu\{w\in g^{-1}(x) : |w|\geq |x|+k\} \leq \nu'\{w\in g^{-1}(x) :
|w|\geq |x|+k\}.$$
Notice that $\langle x\rangle = g^{-1}(x)$, thus $\nu(\langle x\rangle) = \nu'(\langle x\rangle)$ for all $\nu,\nu'\in f^{-1}(\delta_x)$.
The ordered set $\mathbb{O}_{\delta_x}(f)$ is thus, isomorphic to
the following
$$\mathbb{O}_{\delta_x}(f)\cong \{ t\in[0,1]^{n-|x|} :
t_1\leq\ldots\leq t_{n-|x|}\}.$$

\medskip

\begin{prop}
The ordered set $O_k = \{ t\in[0,1]^{k} : t_1\leq\ldots\leq t_k\}$
is an irreducible ordered set.
\end{prop}

Proof: We proceed by induction on $k$. For $k=0$, $O_0$ is a
singleton (by convention, if desired) and for $k=1$ we have that
$O_1 = [0,1]$ is linearly ordered, so we suppose that $k\geq 2$
and that we have an order-isomorphism $\phi:O_k\To P\times Q$. We
denote by the symbols 0 and 1 the minimum and the maximum
respectively of any of the ordered sets $O_k$, $P$ and $Q$ (all
must exist so that $\phi(0)=(0,0)$ and $\phi(1)=(1,1)$). Let
$$\Lambda = \{ t\in O_k : t_2 = t_3=\cdots=t_k=1\} = \{t\in O_k :
\{s : s\geq t\}\text{ is linearly ordered }\}$$

Every element of $\phi(\Lambda)$ must be either of the form
$(x,1)$ or $(1,x)$, since otherwise $\{s :s\geq \phi(\lambda)\}$
cannot be linearly ordered. Moreover, since $\Lambda$ is linearly
ordered, it follows that either $\phi(\Lambda)\subset
P\times\{1\}$ or $\phi(\Lambda)\subset \{1\}\times Q$. We suppose
that $\phi(\Lambda)\subset P\times\{1\}$. Now call $\lambda =
(0,1,\ldots,1)\in \Lambda$ and $\phi(\lambda) = (u,1)$. We have
that
$$O_{k-1}\cong \phi\{t\in O_k :
t\leq\lambda\} = \{s\in P : s\leq u\}\times Q,$$
so by the inductive hypothesis, either $|Q|=1$ (which would finish
the proof) or $u=0$. So we suppose that $u=0$, which implies that
$Q\cong O_{k-1}$ and also that $\phi(\Lambda) = P\times\{1\}$
(because we found that $\phi(\lambda)=(0,1)\in\phi(\Lambda)$ and
this is an upwards closed set). Thus $Q\cong O_{k-1}$ and
$P\cong\Lambda\cong [0,1]$, and it remains to show that
$O_k\not\cong O_{k-1}\times [0,1]$. The reason is that the
elements $p=((0,1,\ldots,1),1)$ and $q=((1,\ldots,1),0)$ are two
incomparable elements of $O_{k-1}\times [0,1]$ with the property
that $\{t: t\geq p\}$ and $\{t: t\geq q\}$ are linearly ordered.
However we noticed that the set $\Lambda$ of points with this
property in $O_k$ is linearly
ordered.$\qed$

Since one of our announced objectives was to show that
$P(\sigma_n(\kappa))$ is not homeomorphic to $P(\sigma_m(\kappa))$
for $n\neq m$ let us make explicit now why this is true. It is
enough to notice that the irreducible ordered sets $O_k = \{t\in
[0,1]^k : t_1\leq\cdots\leq t_k\}$ which appear in the fiber
orders of these spaces are not order-isomorphic for different
values of $k$, since for $n<m$, $O_m$ does not appear as the fiber
order of any point of $P(\sigma_n(\kappa))$. This can be realized
in many different ways. We propose to the reader one of them.
Consider
$$e=(0,0,\ldots,1) = \max\{t\in O_k : \{s : s\leq t\}\text{ is
linearly ordered}\}.$$ Then $O_{k-1}\cong\{t\in O_k : t_k=1\} =
\{t: t\geq e\}$. This argument shows that $O_{k-1}$ can be
obtained in an intrinsic way from $O_k$, and thus if $O_k\cong
O_j$, then $O_{k-1}\cong O_{j-1}$. The inductive
repetition of such argument leads to contradiction if $k\neq j$.

Finally, let us show the appropriate parts of Theorems~\ref{powerBall} and~\ref{productball}.
We get easily that $\sigma_n(\kappa)$ satisfies the assumptions of Proposition~\ref{assumdecomp}(i) with $O=O_n$ and $x=\emptyset$. Further, $A(\kappa)=\sigma_1(\kappa)$ satisfies the assumptions of Proposition~\ref{assumdecomp}(ii) with $O=O_1=[0,1]$ and $x=\emptyset$.

\subsection{The spaces $P([0,\omega_1]^n)$ and
$P(A(\kappa)^n)$}\label{sectPAkn}
In this section we shall prove Theorem~\ref{polyhedric} and the appropriate part of Theorems~\ref{powerBall} and~\ref{productball}.

The fiber orders of the two spaces of probability measures from
the title can be computed in the same way. For $M\supset N$, let
$p_{MN}:A(M)\To A(N)$ be the continuous surjection given by
$p_{MN}(x)=x$ if $x\in A(N)$ and $p_{MN}(x)=\infty$ otherwise. The
$\sigma$-typical surjection of $P(A(\kappa)^n)$ is of the form
$P(p_{MN}^n):P(A(M)^n)\To P(A(N)^n)$ for $M\supset N$ infinite
countable subets of $\kappa$
such that $M\setminus N$ is infinite.

On the other hand, for countable ordinals $\alpha<\beta$ let
$q_{\beta\alpha}:[0,\beta]\To [0,\alpha]$ be the continuous
surjection given by $q_{\beta\alpha}(\gamma)=\gamma$ for
$\gamma\leq\alpha$, and $q_{\beta_\alpha}(\gamma)=\alpha$ for
$\gamma>\alpha$. The $\sigma$-typical surjection of
$P([0,\omega_1]^n)$ is of the form $P(q_{\beta\alpha}^n):
P([0,\beta]^n)\To P([0,\alpha]^n)$ where $\alpha<\beta$ are
countable limit ordinals. From the point of view of fiber orders
both surjections $p_{MN}$ and $q_{\beta\alpha}$ can be treated
simultaneously since both can be viewed as a surjection $g:K\To L$
satisfying the following properties:

$(\star)$ There exist a point $\varpi\in L$ and a point $m\in
g^{-1}(\varpi)$ such that $|g^{-1}(x)|=1$ for every $x\in
L\setminus\{\varpi\}$, and with respect to the fiber order of
$g^{-1}(\varpi)$, we have that $m<t$ and $t\sim s$ for every
$t,s\in
g^{-1}(\varpi)\setminus\{m\}$.

In the case of $p_{MN}$ we should take $\varpi=\infty$ and
$m=\infty$, while for $q_{\beta\alpha}$, $\varpi=\alpha$ and
$m=\alpha$.

From now on, we shall
concentrate in computing the fiber orders of $P(g^n)$ where
$g:K\To L$ is a continuous surjection satisfying $(\star)$, and
with this information the computation of the fiber orders of
$P(A(\kappa)^n)$ and $P([0,\omega_1]^n)$ will follow
immediately.

We fix $x=(x_1\ldots,x_n)\in L^n$, and we call $R(x) =
\{i\in\{1,\ldots,n\} : x_i = \varpi\}$. First step is to
understand which are the sets of the form $\langle
y^{(1)},\ldots,y^{(k)} \rangle$ in $(g^n)^{-1}(x)$. For every
$y\in (g^n)^{-1}(x)$ we call $S(y) = \{i\in R(x) : y_i > m\}$.

Claim A: $\langle y^{(1)},\ldots, y^{(k)}\rangle = \{z\in
(g^n)^{-1}(x) : \exists j\in\{1,\dots,k\} : S(z) \supset
S(y^{(j)})\}$.

Proof of Claim A: Suppose first that $S(z)\supset S(y^{(j)})$ for
some $j$. Then it follows immediately that $y^{(j)}\leq z$ since
the inequality holds coordinatewise. Thus $z\in \langle
y^{(1)},\ldots, y^{(k)}\rangle$.

Now, for the converse inclusion suppose that for every $j$ there
exists a coordinate $i(j)\in S(y^{(j)})\setminus S(z)$. So
$z_{i(j)}=m$ and $y^{(j)}_{i(j)}>m$. Since all the elements of
$\mathbb{F}_x(g)$ which are greater than $m$ are equivalent, for
every $y^{(j)}_i>m$ we can easily find a neighborhood $W^j_i$ such
that $W = \bigcup_{i,j} g(W^j_i)$ does not contain the image of
any neighborhood of $m$. For every $j\in\{1,\ldots,k\}$ set $U_j =
\{y\in K^n : y_{i(j)}\in W^j_{i(j)}\}$ which is a neighborhood of
$y^{(j)}$. We claim that $g^n(U_1)\cup\cdots\cup g^n(U_k)$
contains no image of a neighborhood of $z$, which will finish the
proof of Claim A. Namely, if $V$ is a neighborhood of $z$ of the
form $V_1\times\cdots\times V_n$ with $V_i$ neighborhood of $z_i$,
then for every $i\in R(x)\setminus S(z)$ we can find a point
$t_i\in g(V_i)\setminus W$. If we take $u\in V$ with $g(u_i)=t_i$
for $i\in R(x)\setminus S(z)$, then
$g^n(u)\not\in g^n(U_1)\cup\cdots\cup g^n(U_k)$.

Claim B: For $\nu,\nu'\in P(g^n)^{-1}(\delta_x)$, we have that
$\nu\leq\nu'$ if and only if for every upwards closed subset of
the power set of $R(x)$, $A\subset 2^{R(x)}$, we have that $\nu\{z
:S(z)\in A\}\leq \nu'\{z : S(z)\in A\}$.

Proof of Claim B: It follows from Claim A that the subsets of
$(g^n)^{-1}(x)$ of the form $\langle y^{(1)},\ldots,
y^{(k)}\rangle$ are exactly the sets of the form $\{z : S(z)\in
A\}$ for some upwards closed family $A$ of subsets of $R(x)$.

As a consequence, for $x\in L^n$ with $|R(x)|=k$, we have that
$\mathbb{O}_{\delta_x}(P(g^n))\cong \{t\in [0,1]^{2^k} :
\sum_{i\in 2^k} t_i = 1\}$ endowed with the order $t\leq s$ if and
only if $\sum_{i\in A}t_i\leq \sum_{i\in A}s_i$ for every upwards
closed subset of $2^k$.

\begin{prop}
Consider the ordered set $P_k =\{t\in [0,1]^{2^k} : \sum_{i\in
2^k} t_i = 1\}$ endowed with the order $t\leq s$ if and only if
$\sum_{i\in A}t_i\leq \sum_{i\in A}s_i$ for every upwards closed
subset of $2^k$. Then $P_k$ is an irreducible ordered set.
\end{prop}

Proof: We proceed by induction on $k$. If $k=1$, then $P_k\cong
[0,1]$. Suppose that we had an isomorphism $\phi:P_k\To Q\times
R$. We shall use the symbols 0 and 1 to denote the minimum and
maximum of any of these ordered sets (notice that that the minimum
of $P_k$ is the characteristic function of the empty set
$0=\chi_{\{\emptyset\}}$, and its maximum is $1=\chi_{\{\{1,\ldots,n\}\}}$).
For every $i\in\{1,\ldots,n\}$ we consider $e^i\in P_k$ the
characteristic function of the singleton $\{\{i\}\}$. Notice that
$\{t\in P_k : t\leq e^i\}$ is linearly ordered since any such $t$
satisfies $\sum_{a\not\subset\{i\}} t_a=0$. Thus $\phi(e^i)$ must
be of the form either $\phi(e^i)=(u_i,0)$ or $\phi(e^i) =
(0,u_i)$. Notice now that
$$\{t\in P_k : t\geq e^i\} = \{t\in P_k
: \sum_{i\in a}t_a = 1\}\cong P_{k-1}$$
and $P_{k-1}$ is irreducible by the inductive hypothesis, so
$u_i=1$ since $\{(r,s) : (r,s)\geq (r_0,s_0)\} = \{r: r\geq
r_0\}\times\{s : s\geq s_0\}$. Hence $\phi(e^i)\in\{(1,0),(0,1)\}$
for every $i\in\{1,\ldots,k\}$. If $k>2$ this is already a
contradiction, so we suppose that $k=2$ and $\phi(e^1) = (1,0)$
and $\phi(e ^2) = (0,1)$. We denote the elements of $P_2$ as
$t=(t_\emptyset,t_{\{1\}},t_{\{2\}},t_{\{1,2\}})$. For every
$\lambda\in [0,1]$, we call $x^\lambda=(1-\lambda,\lambda,0,0)$
and $y^\lambda=(1-\lambda,0,\lambda,0)$ in $P_2$. We have
$x^\lambda\le e^1$ and $y^\lambda\le e^2$ and so $\phi(x^\lambda)
= (r^\lambda,0)$ and $\phi(y^\lambda) = (0,s^\lambda)$ for
suitable $r^\lambda$ and $s^\lambda$. We consider the specific
elements $u=(0.7,0.1,0.1,0.1)$ and $u'=(0.8,0,0,0.2)$ of $P_2$.
Say that $\phi(u)=(r,s)$ and $\phi(u') = (r',s')$. On the one hand
$x^{0.2}$ and $y^{0.2}$ are lower than $u$ and $u'$ so
$r^{0.2}\leq r$, $r^{0.2}\leq r'$ , $s^{0.2}\leq s$ and
$s^{0.2}\leq s'$. On the other hand, if $\lambda>0.2$ then
$x^\lambda\not\leq u$, $x^\lambda\not\leq u'$, $y^\lambda\not\leq
u$, neither $y^\lambda\not\leq u'$. Hence indeed $r=r'=r^{0.2}$
and $s=s'=s^{0.2}$. Thus $\phi(u)=\phi(u')$, a
contradiction.$\qed$

\medskip

Notice that  $P_k$ is not order-isomorphic to $P_{k'}$ for $k\neq
k'$, since the set
\begin{multline*}
H = \{t\in P_k : \{s : s\leq t\}\text{ is linearly ordered}\} \\ =
\{t\in P_k: \exists i\in\{1,\ldots,k\} : t\leq e^i\}
\end{multline*} contains
exactly $k$ many maximal elements: $\{e^1,\ldots,e^k\}$, where
again $e^i\in P_k$ denotes the characteristic function of the
singleton $\{i\}$. This also shows that these irreducible ordered
sets are not isomorphic to the irreducible ordered sets $O_k$
which appeared in the fiber orders of the spaces
$P(\sigma_n(\kappa))$ (for $n>1$), because in those cases the set
of all elements $t$ such that $\{s :s\leq t\}$ is linearly ordered
was a linearly ordered set with precisely one maximal element.

\medskip

The above calculation proves Theorem~\ref{polyhedric}. Further,
both $A(\kappa)^n$ and $[0,\omega_1]^n$ satisfy the assumptions of
Proposition~\ref{assumdecomp}(i) with $O=P_n$ and
$x=(\infty,\dots,\infty)$ resp.
$x=(\omega_1,\dots,\omega_1)$.  This proves the appropriate part
of Theorem~\ref{powerBall}. Finally, $[0,\omega_1]$ satisfies the
assumptions of Proposition~\ref{assumdecomp}(ii) with $O=[0,1]$
and $x=\alpha$ (using the above notation). This proves the appropriate part of
Theorem~\ref{productball}.

\medskip

We have not mentioned it so far but, despite the fact that the
picture of fiber orders is similar, the spaces $P(A(\kappa)^n)$
and $P([0,\omega_1]^n)$ are very different, by other well known
reasons. Namely, $P(A(\kappa)^n)$ is an Eberlein compact, and
hence Fréchet-Urysohn space, so it cannot contain any copy of
$[0,\omega_1]$.

\section{Higher weights}

So far we used the version of spectral theorem that we stated as
Theorem~\ref{ourspectraltheorem} but there is the possibility to
use other versions. For example, for a regular cardinal $\tau$, we
consider $\mathcal{Q}_\tau(K)$ the family of quotients of weight
strictly less than $\tau$, and we call a $\tau$-semilattice to a
subset $\mathcal{S}\subset\mathcal{Q}(K)$ such that the supremum
of every subset of $\mathcal{S}$ of cardinality less than $\tau$
belongs to $\mathcal{S}$. The set $\mathcal{S}$ is cofinal in
$\mathcal{Q}_\tau(K)$ if for every $L\in\mathcal{Q}_\tau(K)$ there
exists $L'\in\mathcal{S}$ with $L\leq L'$. We assume that $\tau$
is a regular cardinal
because otherwise there exists no cofinal $\tau$-semilattice in $\mathcal{Q}_\tau(K)$.

\begin{thm}
Let $K$ be a compact space with weight at least $\tau$. The
intersection of two cofinal $\tau$-semilattices in
$\mathcal{Q}_\tau(K)$ is
a further cofinal $\tau$-semilattice in $\mathcal{Q}_\tau(K)$.
\end{thm}

Proof: It is completely analogous to the proof of
Theorem~\ref{ourspectraltheorem}.$\qed$

\begin{lem}
Let $K$ be a compact space of weight at least $\tau$ and
$\mathcal{S}$ a $\tau$-semilattice in $\mathcal{Q}_\tau(K)$. Then,
$\mathcal{S}$
is cofinal in $\mathcal{Q}_\tau(K)$ if and only if $\sup\mathcal{S} = K$.
\end{lem}

Proof: Suppose $\sup\mathcal{S} = K$. By the same argument as in
the proof of Lemma~\ref{spectralfactor}, every real-valued
continuous function $f\in C(K)$ factors through an element of
$\mathcal{S}$. Now, if $p:K\To L$ is an arbitrary element of
$\mathcal{Q}_\tau(K)$, then we can take an embedding
$L\subset\mathbb{R}^\gamma$ for a cardinal $\gamma<\tau$ and
consider the functions $e_i p:K\To\mathbb{R}$ obtained by
composing with the coordinate functions
$e_i:\mathbb{R}^\gamma\To\mathbb{R}$, $i<\gamma$. For every
$i<\gamma$ we know that there exists $L_i\in\mathcal{S}$ such that
$e_i p$ factors through $L_i$. Finally, this implies that $p$
factors through $L_\infty = \sup\{L_i : i<\gamma\}$, so $L\leq
L_\infty\in\mathcal{S}$.$\qed$

In a similar way as we did with $\sigma$-semilattices spectra, we
can say that the $\tau$-typical surjection of $K$ has a property
$\mathcal P$ if there is cofinal $\tau$-semilattice in which all
the natural surjections have property $\mathcal P$, and when this
happens such a $\tau$-semilattice can be found as a subsemilattice
of any given one. Also similarly, we can talk in this context of
$\mathbb{F}_x^\tau(K)$ and $\mathbb{O}_x^\tau(K)$.

An application can be found in the study of of the space
$P([0,\tau]^n)$ and their finite powers, for $\tau>\omega_1$ a
regular cardinal. The fiber orders of $K=P([0,\alpha]^n)$ for any
ordinal $\alpha\geq\omega_1$ can be computed using very similar
arguments as in Section~\ref{sectPAkn}, and indeed
$\mathbb{O}_{\delta_x}(K)\cong P_k$ where $k$ is the number of
coordinates of $x\in[0,\alpha]^n$ with uncountable cofinality.
Therefore, for $\alpha\geq\omega_1+\omega_1$, $P([0,\alpha]^n)$
does not satisfy the assumptions of Proposition~\ref{assumdecomp}. Indeed, $[0,\alpha]$ has at least two non-$G_\delta$-points and hence $[0,\alpha]^n$ contains several points $x$ with $\mathbb{O}_{\delta_x}(P(K))\cong P_n$. Moreover, $P([0,\alpha]^n)$ does not satisfy even the assumptions of Theorems~\ref{Knodecomposition} and~\ref{Knoroots} -- for $n=1$ it is witnessed by the fact that $\mathbb{O}_{\frac12(\delta_x+\delta_y)}(P([0,\alpha]))=[0,1]^2$ whenever $x$ and $y$ are two distinct points with uncountable cofinality.

However, still it is possible to get decomposition results about
spaces $P([0,\tau]^n)$ using $\tau$-semilattices, since analogues
of Theorems~\ref{Knodecomposition} and \ref{Knoroots} for the
$\tau$-typical surjection hold, with identical proof. There is a
natural cofinal $\tau$-semilattice for $[0,\tau]$: For
$\alpha<\beta$ consider the continuous surjection
$p_{\beta\alpha}:[0,\beta]\To [0,\alpha]$ given by
$p_{\beta\alpha}(x)=x$ for $x\leq\alpha$, and
$p_{\beta\alpha}(x)=\alpha$ for $x>\alpha$. The $\tau$-semilattice
consists of all quotients given by $p_{\tau\alpha}$,
$\alpha<\tau$, and the $\tau$-typical surjection is of the form
$p_{\beta\alpha}$, $\alpha<\beta<\tau$. Thus, the situation is
completely analogous to that of $P([0,\omega_1]^n)$, and we have
the following result:

\begin{thm}
Let $\tau$ be a regular cardinal, $K=P([0,\tau]^n)$,
$x=(x_1,\ldots,x_r)\in [0,\tau]^n$ and $k=|\{i: x_i=\tau\}|$. Then
$\mathbb{O}^\tau_{\delta_x}(K)\cong P_k =\{t\in [0,1]^{2^k} :
\sum_{i\in 2^k} t_i = 1\}$ endowed with the order $t\leq s$ if and
only if $\sum_{i\in A}t_i\leq \sum_{i\in A}s_i$ for every upwards
closed subset of $2^k$.
\end{thm}

Now we get easily the remaing part of Theorem~\ref{powerBall} and
Theorem~\ref{producttau}.

We mention that, answering a question posed to us by R. Deville and G. Godefroy, the ideas of this section are used in \cite{Avilesnumber} to show that there exist $2^\kappa$ many nonhomeomorphic weakly compact convex subsets in $\ell_2(\kappa)$.

\section{Final Remarks and Open Problems}

\begin{question} Let $M(K)$ denote the space of regular Borel measures of variation at most $1$
(that is, the dual ball of the space of continuous functions
$C(K)$) in its weak$^\ast$ topology. We show in this paper that
$M(A(\kappa))$ is not homeomorphic to $P(A(\kappa))$ using fiber
orders. We did not make a systematic study of the fiber orders of
$M(K)$ and this may be interesting. Analysing the relatively easy
case of $A(\kappa)$ it seems that the fiber orders of $M(A(\kappa))$ look
similar to those of $P(A(\kappa))^2$, so we may ask: Is $M(K)\approx P(K)^2$ for each compact space? However, this question has negative answer. Let $K$ be the well-known ``double arrow space''. Then $P(K)$ is first countable (for example by \cite[Proposition 7]{phd}) while $M(K)$ is not first countable as $K$ is not metrizable. Therefore $M(K)\not\approx P(K)^2$. But we still can ask: Is $M(K)\approx P(K)^2$ for compact spaces considered in this paper ($A(\kappa)$, $\sigma_n(\kappa)$ etc.)?
\end{question}
\begin{question} The analysis of the generic fibers
of $B(\kappa)$ yields the same result as for $P(A(\kappa))$,
namely all the non $G_\delta$ points have generic fibers
order-isomorphic to an interval. Are
the spaces $B(\kappa)$ and $P(A(\kappa))$ homeomorphic? In relation with this, it follows from \cite{Kalendaproducts} that $P(A(\kappa))$ is homeomorphic to $P(A(\kappa))\times [0,1]$. Is $B(\kappa)$ homeomorphic to $B(\kappa)\times [0,1]$ or even to any product of two nontrivial spaces?
\end{question}
\begin{question}  In the various spaces of probability measures that we
studied, fiber orders allow us to identify different types of
points. Is this a complete classification? That is, we ask to
determine exactly for which points $x,y\in P(K)$ there exists a
homeomorphism $f:P(K)\To P(K)$ such that $f(x)=y$.
\end{question}
\begin{question} Fiber orders are a good tool to determine whether two spaces are homeomorphic but they do not seem
to help in determining wheter a given space is the continuous
image of another. In \cite{Avilesproducts} the case of the spaces
$B(\kappa)^n$ is studied, but the situation is not clear for the
other spaces studied here. For instance we do not know whether
$P(\sigma_n(\kappa))$ maps onto $P(\sigma_m(\kappa))$ for $n<m$,
and so on. This is related also to the problem of the
$A(\kappa)^\omega$-images, initiated by Benyamini, Rudin and Wage
\cite{BenRudWag} and studied specially by Bell
in~\cite{BellRamsey} and \cite{Bellctightness}. It is proven
in~\cite{AvilesHilbertBall} that $P(A(\kappa))$ and $B(\kappa)$
are continuous images of $A(\kappa)^\omega$, but it is unclear to
us whether $P(\sigma_n(\kappa))$ or $P(A(\kappa)^n)$ are
continuous images of $A(\kappa)^\omega$ for $n>1$.
\end{question}

\section*{Acknowledgement}

We want to express our gratitude to Stevo Todor\v{c}evi\'c for his
valuable suggestions, and also to Richard Haydon, who proposed to
us a shorter proof of Lemma~\ref{numericalimage}. Thanks are also due to Petr Holick\'y, Ji\v{r}\'{\i} Spurn\'y and Miroslav Zelen\'y for helpful comments.

\bibliographystyle{amsplain}

\end{document}